\documentclass[12pt]{article}

\usepackage[a4paper,margin=27mm]{geometry}
\usepackage[T1]{fontenc}
\usepackage{lmodern}
\usepackage{amsmath,amssymb,amsthm,mathtools,mathrsfs}
\usepackage{microtype}
\usepackage{enumitem}
\usepackage{booktabs}
\usepackage{indentfirst}
\usepackage{xcolor}
\usepackage{authblk}
\usepackage{hyperref}
\usepackage{titlesec}

\hypersetup{hidelinks}

\setlist{nosep,leftmargin=2em}

\titleformat{\section}{\large\bfseries}{\thesection}{0.7em}{}
\titleformat{\subsection}{\normalsize\bfseries}{\thesubsection}{0.7em}{}

\newtheorem{theorem}{Theorem}[section]
\newtheorem{proposition}[theorem]{Proposition}
\newtheorem{lemma}[theorem]{Lemma}

\theoremstyle{definition}

\theoremstyle{remark}
\newtheorem{remark}[theorem]{Remark}

\newcommand{\Ct}{\mathcal C_\theta}
\newcommand{\sig}{\sigma}
\newcommand{\W}{W}
\newcommand{\bsig}{\bar\sigma}
\newcommand{\cR}{\mathbb R}
\newcommand{\cE}{\mathcal E}
\newcommand{\cM}{\mathcal M}
\newcommand{\cT}{\mathcal T}
\newcommand{\cL}{\mathcal L}
\newcommand{\dd}{\,d}

\newcommand{\od}{\mathbin{\odot}}

\title{The Capillary Orlicz--Christoffel--Minkowski Flows}
\author[1]{Dongjing Ma\thanks{E-mail: \href{mailto:3261019581@qq.com}{3261019581@qq.com}}}
\author[1]{Peihe Wang\thanks{E-mail: \href{mailto:peihewang@hotmail.com}{peihewang@hotmail.com}}}
\author[2]{Hua Zhu\thanks{E-mail: \href{mailto:zhuhmaths@mail.ustc.edu.cn}{zhuhmaths@mail.ustc.edu.cn}}}

\affil[1]{School of Mathematical Sciences, Qufu Normal University, Qufu 273165, Shandong Province, People's Republic of China}
\affil[2]{School of Mathematics and Physics, Southwest University of Science and Technology, Mianyang 621010, Sichuan Province, People's Republic of China}

\date{}
\begin{document}
\maketitle

\begin{abstract}
We study a capillary Orlicz--Christoffel--Minkowski flow on a non-obtuse spherical cap for
$1\le k\le n-1$:
\[
\partial_t h=
\frac{h^2\varphi(\ell/h)}{f}\,
\bsig_k(\nabla^2h+h\sig)-h,
\qquad
\nabla_\mu h=\cot\theta\,h.
\]
In the supercritical regime
\[
a_\varphi(s)=\frac{s\varphi'(s)}{\varphi(s)}
\ge k+1+\delta,
\qquad
s a_\varphi'(s)\ge0,
\]
we obtain time-uniform $C^0$ and $C^1$ estimates and a closed scalar equation for the speed ratio.
A dynamic Orlicz square-completion identity, combined with inverse-concavity, yields a viscosity
estimate for the minimum principal radius. Under a strict dynamic Orlicz-convexity condition and
\[
\frac{f_\mu}{f}\ge\frac{n-k}{n}\cot\theta,
\]
we prove time-uniform two-sided principal-radius bounds on every smooth
existence interval for $0<\theta\le\pi/2$ and deduce long-time
existence and smooth convergence to a stationary solution. The boundary
regularity step uses a quantitative interior concave estimate through
Neumann reflection and a concave uniformly elliptic extension.
Under the weaker condition
$(f\ell^{k-1})_\mu\ge0$, strict convexity is preserved on every smooth existence interval.
For $k=1$, the transformed scalar equation is uniformly parabolic for every
$0<\theta<\pi$. Classical parabolic estimates then give an all-angle
scalar flow theorem, without asserting geometric convexity preservation
in the obtuse range.
\end{abstract}


\section{Introduction}

Prescribed area-measure and curvature equations form a central part of the Brunn--Minkowski
theory and fully nonlinear geometric analysis. In the closed setting, the classical
Christoffel--Minkowski problem is naturally expressed through elementary symmetric functions
of the principal radii; its fully nonlinear analytic theory is closely related to Hessian
equations and the convexity principle \cite{CNS,GuanMa,GuanMaZhou}. The $L_p$ extension is
linked to the Brunn--Minkowski--Firey theory \cite{Lutwak93}. The Orlicz--Minkowski problem
and the broader Orlicz--Brunn--Minkowski theory were developed in the works
\cite{GardnerHugWeil,HaberlLutwakYangZhang,XiJinLeng}. Curvature flows provide an
effective approach to these problems: Ivaki \cite{Ivaki2019} studied flows driven by
$\sigma_k$ of the principal radii and the support function, while Bryan--Ivaki--Scheuer
developed broad parabolic approaches to curvature equations and Orlicz--Minkowski flows
\cite{BISOrlicz,BISParabolic}. Their inverse-concavity and viscosity methods provide an
important model for the present work \cite{BISConstantRank}.

For capillary hypersurfaces in the Euclidean half-space, the support function is defined on a
spherical cap and satisfies a Robin boundary condition. The capillary Minkowski problem was
introduced and solved in \cite{MWWMinkowski}, with an $L_p$ extension developed in
\cite{MWWLpMinkowski}. The intermediate-order capillary Christoffel--Minkowski problem has
been studied from complementary viewpoints in \cite{HISCapillaryCM,MWWChristoffel}, and
related capillary $L_p$ curvature and Christoffel--Minkowski problems were treated in
\cite{HuIvaki,HuIvakiCurvature,MWWLpCurvature}. Capillary Alexandrov--Fenchel and
Minkowski-type inequalities have also been obtained by curvature-flow methods
\cite{WWXAF,WWXMinkowski}. At the Monge--Amp\`ere endpoint $k=n$, the capillary
Orlicz--Minkowski problem was developed in \cite{WangZhu}, while Li and Liu
\cite{LiLiuFlow} studied the capillary Orlicz--Minkowski flow in the upper half-space.

The present paper concerns the intermediate orders $1\le k\le n-1$. In this sense, it is
complementary to the determinant-level flow of Li--Liu and is a dynamic companion to the
elliptic capillary Orlicz--Christoffel--Minkowski problem studied by Zhu
\cite{ZhuStatic}. Let
\[
\cR^{n+1}_+=\{x\in\cR^{n+1}:x_{n+1}>0\},
\qquad
e=-E_{n+1},
\]
and let $\Ct$ be the unit capillary spherical cap with non-obtuse contact angle
\[
0<\theta\le\frac{\pi}{2}.
\]
Its support function is
\[
\ell(\xi)=\sin^2\theta+\cos\theta\langle\xi,e\rangle,
\]
and
\begin{equation}\label{ell-identities}
\nabla^2\ell+\ell\sig=\sig,
\qquad
\nabla_\mu\ell=\cot\theta\,\ell.
\end{equation}
For $0<\theta\le\pi/2$ one has
\begin{equation}\label{ell-bounds}
0<1-\cos\theta\le \ell\le\sin^2\theta,
\end{equation}
with $\ell\equiv1$ when $\theta=\pi/2$.
For a positive capillary support function $h$ we put
\[
\W[h]=\nabla^2h+h\sig.
\]
The eigenvalues of $\W[h]$ are the principal radii. We use the normalization
\[
\bsig_k(A)=\binom nk^{-1}\sigma_k(A).
\]
The stationary Orlicz--Christoffel--Minkowski equation is
\begin{equation}\label{stationary-eqn}
h\,\varphi(\ell/h)\,\bsig_k(\W[h])=f
\quad\hbox{in }\Ct,
\qquad
\nabla_\mu h=\cot\theta\,h
\quad\hbox{on }\partial\Ct.
\end{equation}

We study the flow without a global normalization term
\begin{equation}\label{main-flow}
\partial_t h=
\frac{h^2\varphi(\ell/h)}{f}\,
\bsig_k(\W[h])-h,
\qquad
\nabla_\mu h=\cot\theta\,h.
\end{equation}
Its stationary points are exactly the solutions of \eqref{stationary-eqn}. In the power case
$\varphi(s)=s^p$, the speed becomes
\[
\frac{\ell^p}{f}h^{2-p}\bsig_k(\W[h])-h,
\]
which is the capillary analogue of the closed $L_p$ Christoffel--Minkowski flows in
\cite{BISOrlicz,Ivaki2019}.

Our first observation is variational. The functional
\[
\cM_k[h]=\frac1{k+1}\int_{\Ct}h\,\bsig_k(\W[h])\dd\sig
\]
has first variation
\[
D\cM_k[h](\eta)=\int_{\Ct}\eta\,\bsig_k(\W[h])\dd\sig
\]
for Robin variations. If
\[
J_\varphi(s)=\int_1^s\frac{\dd\tau}{\tau\varphi(1/\tau)},
\]
then
\[
\cE_{f,\varphi}[h]
=
\cM_k[h]-\int_{\Ct}fJ_\varphi(h/\ell)\dd\sig
\]
has Euler--Lagrange equation \eqref{stationary-eqn}, and \eqref{main-flow} is a
gradient ascent flow for this functional with positive mobility.

The central scalar quantity is
\begin{equation*}
Q=
\frac{h\varphi(\ell/h)\bsig_k(\W[h])}{f}.
\end{equation*}
We prove that $Q$ satisfies a closed scalar parabolic equation with a homogeneous Neumann
boundary condition. In the supercritical range
\begin{equation}\label{supercritical}
a_\varphi(s):=\frac{s\varphi'(s)}{\varphi(s)}
\ge k+1+\delta
\quad(s>0)
\end{equation}
this yields time-uniform upper and lower bounds for $Q$ and, eventually, exponential convergence
$Q\to1$.

The main new difficulty is the principal-radius estimate. For the dynamic problem it is useful to
introduce
\[
H=f\ell^{k-1},\qquad q_0=\log\ell,\qquad K=\log H.
\]
For $s>0$ put
\[
b(s)=a_\varphi(s)-2,\qquad
d(s)=s a_\varphi'(s),
\]
and
\begin{equation*}
A_\varphi^{\rm fl}(s)
=
b(s)+d(s)+\frac{b(s)^2}{k+1}.
\end{equation*}
The dynamic structural tensor is
\begin{equation}\label{Tflow-def}
\begin{aligned}
\cT_{\varphi,f}(s)
:={}&
-\nabla^2K-2\,dq_0\od dK
+\frac{b(s)+d(s)}
{(k+1)A_\varphi^{\rm fl}(s)}\,dK\otimes dK\\
&+
\frac{b(s)+k+1}{\ell}\,\sig .
\end{aligned}
\end{equation}
Here $\alpha\od\beta=\frac12(\alpha\otimes\beta+\beta\otimes\alpha)$. In the power case
$\varphi(s)=s^p$, \eqref{Tflow-def} reduces to
\[
\cT_{\varphi,f}
=
(p+k-1)
\frac{
\nabla^2(\ell H^{-1/(p+k-1)})
+
\ell H^{-1/(p+k-1)}\sig
}{
\ell H^{-1/(p+k-1)}
},
\]
recovering the spherical convexity exponent in \cite{Ivaki2019}.

We now state the quantitative theorem. The boundary assumption is not
claimed to be optimal. The curvature estimate is an a priori statement
on a smooth existence interval. Its passage to global smooth existence
uses the classical parabolic estimates proved in Lemma~\ref{lem-parabolic}
and the quantitative interior concave estimate in Lemma~\ref{lem-EKH}.
The reflection argument below reduces the boundary regularity step to
the latter estimate.

\begin{theorem}[Radius bounds and global flow]\label{thm-main}
Let
\[
0<\theta\le\frac{\pi}{2},
\qquad
1\le k\le n-1.
\]
Let $f\in C^\infty(\overline\Ct)$ and
$\varphi\in C^\infty((0,\infty))$ be positive. Assume \eqref{supercritical},
\begin{equation}\label{aphi-monotone}
s a_\varphi'(s)\ge0,
\end{equation}
and suppose that there is $\tau_0>0$ such that
\begin{equation}\label{dynamic-convexity}
\cT_{\varphi,f}(s)\ge\tau_0\sig
\qquad\hbox{for every }s>0.
\end{equation}
Assume also
\begin{equation}\label{boundary-quant}
\frac{f_\mu}{f}
\ge
\frac{n-k}{n}\cot\theta
\qquad\hbox{on }\partial\Ct.
\end{equation}
Let $h_0$ be a smooth positive strictly convex capillary support function
satisfying the compatibility hierarchy specified in
Subsection~\ref{subsec-local}. Every smooth positive solution of
\eqref{main-flow} issuing from $h_0$ remains strictly convex and satisfies
\[
0<c\sig\le \W[h(\cdot,t)]\le C\sig
\qquad(0\le t<T),
\]
on its smooth existence interval $[0,T)$, with $c,C$ independent of $T$.
Moreover, \eqref{main-flow} admits a unique smooth
solution for all $t\ge0$, and there exists a smooth strictly convex
capillary support function $h_\infty$ such that
\[
h(\cdot,t)\longrightarrow h_\infty
\quad\hbox{in }C^\infty(\overline\Ct)
\quad\hbox{as }t\to\infty.
\]
The convergence of $Q$ to $1$ is exponential, and $h_\infty$ solves
\eqref{stationary-eqn}.
\end{theorem}

The natural capillary boundary condition arising in the static problem is the weaker condition
\begin{equation}\label{Hmu-natural}
H_\mu\ge0,
\qquad H=f\ell^{k-1},
\end{equation}
or equivalently
\[
\frac{f_\mu}{f}+(k-1)\cot\theta\ge0.
\]
Although the argument below does not yield a time-uniform quantitative non-collapsing estimate
from \eqref{Hmu-natural} alone, strict convexity is preserved.

\begin{theorem}[Finite-time preservation under the natural boundary sign]\label{thm-preserve}
Let $1\le k\le n-1$ and $0<\theta\le\pi/2$, and let $f$ and
$\varphi$ be smooth and positive as above. Assume
\eqref{supercritical}, \eqref{aphi-monotone},
$\cT_{\varphi,f}(s)\ge0$ for every $s>0$, and
\eqref{Hmu-natural}. If a smooth positive solution of
\eqref{main-flow} exists on $[0,T]$ and
$\W[h_0]>0$, then
\[
\W[h(\cdot,t)]>0
\qquad\hbox{for every }0\le t\le T.
\]
\end{theorem}

\begin{theorem}[Full-angle scalar flow for $k=1$]\label{thm-k1-full-angle}
Let $0<\theta<\pi$, let $f\in C^\infty(\overline\Ct)$ and
$\varphi\in C^\infty((0,\infty))$ be positive, and assume
\[
a_\varphi(s)\ge2+\delta
\qquad (s>0)
\]
for some $\delta>0$. Let $h_0>0$ be smooth and satisfy the Robin compatibility conditions, and
assume
\[
\bsig_1(W[h_0])>0.
\]
Then the scalar flow
\[
h_t=\frac{h^2\varphi(\ell/h)}{f}\,\bsig_1(W[h])-h,
\qquad h_\mu=\cot\theta\,h,
\]
has a unique smooth positive solution for all $t\ge0$. Its speed ratio satisfies
\[
0<c_Q\le Q(\xi,t)\le C_Q,
\qquad
\|Q(\cdot,t)-1\|_{L^\infty}\le Ce^{-ct},
\]
and there is a unique positive stationary solution $h_\infty$ with
$\bsig_1(W[h_\infty])>0$ such that
\[
h(\cdot,t)\longrightarrow h_\infty
\quad\hbox{in }C^\infty(\overline\Ct).
\]
For $\theta>\pi/2$ this theorem is a scalar parabolic result; it does not assert that
$W[h(\cdot,t)]$ remains positive definite.
\end{theorem}

\begin{remark}[Contact-angle range]\label{rem-angle}
The support-function identities \eqref{ell-identities}, the Robin first variation, the $C^0$ and
$C^1$ estimates, the $Q$-equation, and the interior Orlicz square-completion are meaningful for
every fixed contact angle $0<\theta<\pi$ because $\ell$ remains positive. The present high-order
global theorem is restricted to the non-obtuse range because a tangential minimum-radius
support $v=hW(e,e)$ satisfies
\[
v_\mu=h\cot\theta\,\lambda_n.
\]
For $0<\theta<\pi/2$ this is positive, while at $\theta=\pi/2$ it is zero; both signs are
compatible with the boundary barrier used in the non-collapsing argument. Thus the right-angle case
requires no new curvature mechanism. For $\pi/2<\theta<\pi$, however, $\cot\theta<0$ and this
boundary sign is reversed. The present argument therefore does not control the obtuse-angle
boundary. This is the geometric obstruction in the present proof to
extending Theorems~\ref{thm-main}--\ref{thm-preserve} to all contact
angles. For $k=1$ the transformed scalar equation is uniformly
parabolic after the support bound, giving the separate full-angle
result above. Preservation of
geometric strict convexity in the obtuse range is not asserted here.
\end{remark}

The paper is organized into five sections. Section~2 develops the capillary support-function
geometry, Newton-tensor identities, short-time flow, variational structure, and the speed-ratio
equation. Section~3 establishes the $C^0$, $C^1$, and $Q$ estimates and derives the dynamic Orlicz
square completion. Section~4 proves preservation of strict convexity and the quantitative
non-collapsing and two-sided principal-radius estimates. Section~5 contains the long-time
existence and convergence arguments, the independent stationary uniqueness
proof, the all-angle scalar result for $k=1$, and the
discussion of the natural boundary condition and the remaining obtuse-angle problem.

\section{Geometric and variational framework}\label{sec-framework}
\subsection{Capillary support functions and Newton tensors}

Consider
\[
\Ct=\{\xi\in\cR^{n+1}_+:|\xi-\cos\theta\,e|=1\},
\qquad e=-E_{n+1},
\]
equipped with its round metric $\sig$. We write $\mu$ for the outward unit conormal of
$\partial\Ct\subset\Ct$. The function $\ell$ in \eqref{ell-identities} is positive for every fixed
$0<\theta<\pi$. Indeed, if $t=\langle\xi,e\rangle$, then
$t\in[\cos\theta-1,0]$ on $\overline\Ct$, and hence
\[
\ell(\xi)=\sin^2\theta+\cos\theta\,t.
\]
Consequently
\begin{equation}\label{ell-full-angle-bounds}
0<\ell_-(\theta):=
\min\{1-\cos\theta,\sin^2\theta\}
\le\ell\le
\max\{1-\cos\theta,\sin^2\theta\}
=:\ell_+(\theta).
\end{equation}
In the non-obtuse range used in the high-order main theorem this reduces to
\eqref{ell-bounds}.

If a strictly convex capillary hypersurface is parametrized by its capillary Gauss map, its position
vector can be reconstructed from the support function by
\[
X(\xi)=\nabla h+h(\xi-\cos\theta\,e).
\]
The capillary angle condition becomes
\begin{equation}\label{Robin-h}
\nabla_\mu h=\cot\theta\,h.
\end{equation}
The radius tensor is
\[
\W[h]=\nabla^2h+h\sig.
\]
Strict convexity is equivalent to $\W[h]>0$.

For a symmetric matrix $A$, let $\sigma_j(A)$ denote the $j$-th elementary symmetric function and set
\[
\bsig_j(A)=\binom nj^{-1}\sigma_j(A).
\]
The G\aa rding cone is
\[
\Gamma_k=\{\lambda:\sigma_j(\lambda)>0,\ 1\le j\le k\}.
\]
We use
\[
S^{ij}(A)
=
\frac{\partial\bsig_k(A)}{\partial A_{ij}}.
\]
When $A=\W[h]$, the Codazzi property implies
\begin{equation}\label{Newton-divfree}
\nabla_iS^{ij}=0.
\end{equation}
Euler homogeneity gives
\begin{equation}\label{Euler-S}
S^{ij}W_{ij}=k\bsig_k(W).
\end{equation}
% R01 | Inverse-concavity terminology: optional clarification, not correction
% Status: verified clarification; original wording is valid under BIS Definition 1.3
The function $\bar\sigma_k^{1/k}$ is concave on $\Gamma_k$ and
inverse-concave in the sense of \cite[Definition~1.3]{BISConstantRank}.
Here the inverse-concavity test is applied to positive definite matrices:
the map $A\mapsto-\bar\sigma_k(A^{-1})^{1/k}$ is concave there.
All uses of the corresponding differential inequality below take place
while $W[h]>0$.

The Robin condition has two useful consequences. In a boundary orthonormal frame
$\{e_\alpha\}_{\alpha<n}$, $e_n=\mu$, one has
\begin{equation}\label{W-mixed}
W_{\alpha n}=0
\end{equation}
and
\begin{equation}\label{W-normal-derivative}
W_{\alpha\alpha;\mu}
=
\cot\theta\,(W_{nn}-W_{\alpha\alpha}).
\end{equation}
% R02 | Self-contained Robin boundary identities
% Status: verified local derivation
For completeness, let $\kappa=\cot\theta$.  The outward conormal
satisfies $\nabla_{e_\alpha}\mu=\kappa e_\alpha$.  Differentiating
$h_\mu=\kappa h$ tangentially therefore gives
\[
0=e_\alpha(h_\mu-\kappa h)
 =\nabla^2h(e_\alpha,\mu),
\]
and hence $W_{\alpha n}=0$.  The tensor $W=\nabla^2h+h\sigma$ is
Codazzi on the unit sphere.  At a boundary point choose the tangential
frame to be geodesic for the induced boundary metric.  Then
$\nabla_{e_\alpha}e_\alpha=-\kappa\mu$ at that point, and
\[
(\nabla_\mu W)(e_\alpha,e_\alpha)
 =(\nabla_{e_\alpha}W)(\mu,e_\alpha)
 =\kappa W_{nn}-\kappa W_{\alpha\alpha}.
\]
This proves both boundary identities, including the case
$\kappa=0$.

It is often useful to divide by $\ell$:
\[
u=\frac h\ell.
\]
Then \eqref{ell-identities} and \eqref{Robin-h} imply
\begin{equation}\label{Neumann-u}
u_\mu=0,
\end{equation}
and
\begin{equation}\label{Btheta}
W[h]
=
\ell\nabla^2u+\nabla\ell\otimes\nabla u
+\nabla u\otimes\nabla\ell+u\sig.
\end{equation}

\subsection{The flow, parabolic estimates, and the variational structure}\label{subsec-local}

We record and prove the parabolic results used in the continuation
arguments, including their precise classical sources.

\begin{lemma}[Parabolic estimates and local existence]\label{lem-parabolic}
On the fixed smooth cap with $0<\theta<\pi$, the following statements
hold. The linear estimates use fixed uniform ellipticity constants.
\begin{enumerate}
\item A bounded smooth solution of a uniformly parabolic linear
nondivergence equation with bounded measurable principal coefficients,
bounded drift, bounded zeroth-order coefficient, and bounded source
has a uniform parabolic H\"older
estimate on smaller positive-time cylinders. The same estimate holds
up to the boundary under homogeneous Neumann data.
\item Let $v$ be a bounded smooth solution of a linear equation
with homogeneous boundary data $Bv=0$. Let $0<\gamma<1$.
If its coefficients and source are uniformly
$C^{\gamma,\gamma/2}$, and the fixed boundary operator is smooth
and uniformly oblique of Neumann or Robin type, then the solution
has uniform $C^{2+\gamma,1+\gamma/2}$ estimates on smaller
positive-time strips, controlled by the source norm and
$\|v\|_\infty$. For inhomogeneous data $Bv=\Phi$, add the boundary
norm $\|\Phi\|_{C^{1+\gamma,(1+\gamma)/2}}$.
Higher-order versions hold with the corresponding coefficient,
source, and boundary-data norms.
\item A smooth fully nonlinear scalar speed whose Hessian derivative
is positive definite near the initial two-jet, together with a fixed
smooth linear uniformly oblique boundary operator $B$, the homogeneous
condition $Bh=0$, and smooth compatible initial data, has a unique
smooth short-time solution. Compatibility means all identities obtained
by differentiating $Bh=0$ in time and using the equation.
A smooth compatible
uniformly parabolic endpoint jet can be used to restart the solution.
\end{enumerate}
The positive-time constants depend only on the indicated data norms,
ellipticity, the fixed geometry and boundary operator, and strip
lengths and separation. They are uniform under time translation.
\end{lemma}
\begin{proof}
For the first assertion, absorb any bounded zeroth-order term $dv$
into the source. Apply Krylov--Safonov
\cite[Theorems 4.2--4.3, p.~164]{KrylovSafonov}; bounded forcing
controls the required $L^{n+1}$ norm on each fixed cylinder.
In Fermi coordinates a Neumann solution satisfies $v_r=v_{yr}=0$
at the boundary, so its even extension is a strong
$W^{2,1}_{n+1}$ solution. Reflect the principal matrix by $RAR$
and the drift by $Rb$, where $R=\operatorname{diag}(1,\ldots,1,-1)$.
The extended coefficients remain bounded and uniformly elliptic;
coefficient jumps are allowed in the cited theorem. Its backward
parabolic distance gives the asserted positive-time estimates.

For the second assertion, use Ladyzhenskaya--Solonnikov--Ural'tseva
\cite[Chapter IV, Theorem 5.3, pp.~320--321]{LSU} for problem (5.4).
Stereographic projection from the omitted south pole identifies the
fixed cap with a bounded Euclidean ball. Since $\theta<\pi$, the map,
the metric coefficients, and their inverses are smooth on its closure.
Uniform parabolicity and obliqueness are preserved with fixed constants,
so the cited Euclidean-cylinder theorem applies.
For $l=\gamma\in(0,1)$, the coefficients of the equation and its source are
in $H^{l,l/2}$, the boundary coefficients and inhomogeneous
first-order boundary data are in $H^{l+1,(l+1)/2}$, and the solution
is in $H^{l+2,l/2+1}$. The fixed cap and boundary operator satisfy
the required smoothness and obliqueness conditions. A smooth time
cutoff $\chi$ vanishing near the earlier end of a strip removes initial
compatibility issues. The cutoff equation has source
$\chi r+\chi'v$ when $Lv=r$, and boundary data $\chi\Phi$.
In the applications below, the lower H\"older norm of $v$ is already
controlled. The estimate in terms of $\|v\|_\infty$ alone follows
by interpolation and nested-strip absorption: on nested strips
$J'\Subset_p J$, with the same terminal time,
\[
 \|v\|_{2+\gamma;J'}\leq
 \varepsilon\|v\|_{2+\gamma;J}
 +C_\varepsilon\operatorname{gap}(J',J)^{-q}
  \bigl(\|v\|_\infty+\|r\|_{\gamma;J}
                    +\|\Phi\|_{1+\gamma;\partial J}\bigr),
\]
where the displayed norms are parabolic H\"older norms, the boundary norm is on
the lateral boundary, and $q<\infty$.
Choose geometrically decreasing gaps and then a fixed
$\varepsilon<2^{-q}$. Iteration gives a convergent geometric series;
its final remainder tends to zero since the solution is smooth on
the fixed outer strip. This proves the asserted bound without a
prior estimate for the highest derivatives of the solution.
Applying the same theorem with larger noninteger
$l$ proves the higher-order statements. The constants depend on the
controlled norms and strip geometry, not on the location of the strips in time.

For local existence, fix $0<\alpha<1$, write the speed as
$\mathcal N(h)$, and set
$h^*=h_0+t\mathcal N(h_0)$. The first two compatibility identities
give $Bh^*=0$. On a sufficiently short interval the linearization
$L_*v=v_t-D\mathcal N(h^*)v$ is uniformly parabolic.
Use the usual unscaled parabolic H\"older norms and the zero-initial spaces
\[
 X_T=\{v\in C^{2+\alpha,1+\alpha/2}:Bv=0,
              \ v(0)=v_t(0)=0\},\qquad
 Y_T=\{f\in C^{\alpha,\alpha/2}:f(0)=0\}.
\]
The source's zero-initial conditions on p.~321 require $f(0)=0$
for $0<l=\alpha<1$, and the zero boundary source satisfies its
trace condition. Conversely the equation with $v(0)=0$ implies
$v_t(0)=f(0)=0$. Thus
\cite[Chapter IV, Theorem 5.4 and (5.14), p.~322]{LSU}
gives an inverse for $L_*:X_T\to Y_T$ whose norm remains bounded
as $T\downarrow0$. The smooth coefficients of $L_*$ are restrictions
of a fixed family on a sufficiently small interval $[0,T_0]$.
The residual $R_0=\mathcal N(h^*)-h_t^*$ belongs to $Y_T$ and
$\|R_0\|_{Y_T}\leq CT^{1-\alpha/2}\to0$: the residual is smooth
and vanishes at $t=0$. The nonlinear remainder
\[
 R(v)=\mathcal N(h^*+v)-\mathcal N(h^*)-D\mathcal N(h^*)v
\]
vanishes initially because the entire spatial two-jet of $v$ vanishes
there, and satisfies
\[
 \|R(v)-R(w)\|_{Y_T}
 \leq C(\|v\|_{X_T}+\|w\|_{X_T})\|v-w\|_{X_T}.
\]
The H\"older product and composition constants in this estimate are
independent of small $T$. Write $C_*$ for the inverse bound and
$C_R$ for the remainder constant. First choose a fixed radius
$\rho$ within the uniformly parabolic neighborhood with
$2C_*C_R\rho<1/2$, and then choose $T$ so small that
$C_*\|R_0\|_{Y_T}\leq\rho/2$.
The map $v\mapsto L_*^{-1}(R_0+R(v))$ preserves this ball and
is a contraction. Its fixed point gives the local solution.
The higher compatibility identities allow the time and tangential
boundary bootstrap, initiated by difference quotients and the linear
theorem, up to the initial time. Uniform ellipticity then recovers
the missing normal derivatives from the equation. This gives
smoothness. The difference of two local solutions satisfies a linear
uniformly parabolic oblique problem with zero initial and boundary
data, so linear uniqueness applies.
The same construction at a smooth uniformly parabolic endpoint
gives continuation.
\end{proof}

We consider
\begin{equation*}
h_t=\Theta\,S-h,
\qquad
\Theta=\frac{h^2\varphi(\ell/h)}f,
\qquad
S=\bsig_k(W[h]).
\end{equation*}
% R03 | Compatibility and short-time applicability
Here and below, put $B=\partial_\mu-\cot\theta$ and let
$H_j=\partial_t^jh|_{t=0}$ be determined recursively by the flow,
starting with $H_0=h_0$. Compatibility means $BH_j=0$ for every
$j\ge0$. The first nontrivial compatibility condition is
\[
\bigl(\Theta_0\bar\sigma_k(W[h_0])-h_0\bigr)_\mu
 =\cot\theta\,
   \bigl(\Theta_0\bar\sigma_k(W[h_0])-h_0\bigr),
\]
or equivalently $(Q_0)_\mu=0$.

We record the conditions for short-time parabolic existence explicitly.
In a sufficiently small neighborhood of the initial two-jet,
$h>0$ and $W[h]>0$.  Compactness and smoothness imply that
\[
F^{ij}:=\frac{\partial}{\partial h_{ij}}
 \left(\Theta\bar\sigma_k(W[h])-h\right)
 =\Theta S^{ij}
\]
has eigenvalues bounded above and below by positive constants.
The boundary operator
$B(x,z,p)=\langle\mu(x),p\rangle-\cot\theta\,z$
is smooth and uniformly oblique, since $B_p\cdot\mu=1$.
For example, at $h_0$ the spatial linearization of the speed is
\[
v\longmapsto
\Theta_0S_0^{ij}v_{ij}
+\left[(2-a_\varphi(\ell/h_0))Q_0
       +\Theta_0\operatorname{tr}(S_0^{ij})-1\right]v,
\]
with boundary condition $v_\mu-\cot\theta\,v=0$.
Thus the local-existence hypotheses of Lemma~\ref{lem-parabolic}
are satisfied, giving smooth local existence and uniqueness for
compatible initial data.
Positivity of $h$ and $W[h]$ persists for a short time by continuity.
At $\theta=\pi/2$ the boundary operator is homogeneous Neumann.

\subsubsection{First variation and monotone energy}

Define
\begin{equation}\label{Mk-def}
\cM_k[h]
=
\frac1{k+1}\int_{\Ct}h\,\bsig_k(W[h])\dd\sig.
\end{equation}

\begin{lemma}[Robin first variation]\label{lem-first-variation}
If $\eta_\mu=\cot\theta\,\eta$ on $\partial\Ct$, then
\[
D\cM_k[h](\eta)
=
\int_{\Ct}\eta\,\bsig_k(W[h])\dd\sig.
\]
\end{lemma}

\begin{proof}
Differentiating \eqref{Mk-def} produces
\[
\frac1{k+1}\int_{\Ct}
\left\{
\eta\,\bsig_k(W)
+
hS^{ij}(\eta_{ij}+\eta\sig_{ij})
\right\}\dd\sig.
\]
Use \eqref{Newton-divfree} twice. The boundary term is
\[
\int_{\partial\Ct}
S^{\mu\mu}(h\eta_\mu-\eta h_\mu)\dd s,
\]
because \eqref{W-mixed} implies that the mixed Newton-tensor components vanish. The common Robin
condition cancels this term. The interior term is
\[
\frac1{k+1}\int_{\Ct}\eta S^{ij}W_{ij}\dd\sig
=
\frac{k}{k+1}\int_{\Ct}\eta\bsig_k(W)\dd\sig,
\]
and the claim follows.
\end{proof}

Let
\[
J_\varphi(s)
=
\int_1^s\frac{\dd\tau}{\tau\varphi(1/\tau)}
\]
and define
\begin{equation*}
\cE_{f,\varphi}[h]
=
\cM_k[h]
-
\int_{\Ct}fJ_\varphi(h/\ell)\dd\sig.
\end{equation*}
Then
\[
D\cE_{f,\varphi}[h](\eta)
=
\int_{\Ct}
\eta\left(
S-\frac{f}{h\varphi(\ell/h)}
\right)\dd\sig.
\]
Since
\[
h_t
=
\Theta
\left(
S-\frac{f}{h\varphi(\ell/h)}
\right),
\]
we have the exact identity
\begin{equation*}
\frac{\dd}{\dd t}\cE_{f,\varphi}[h(t)]
=
\int_{\Ct}
\frac{f}{h^2\varphi(\ell/h)}\,h_t^2\dd\sig
\ge0.
\end{equation*}

\subsubsection{The speed ratio}

Set
\begin{equation}\label{Q-def}
Q=\frac{\Theta S}{h}
=
\frac{h\varphi(\ell/h)S}{f}.
\end{equation}
Then
\begin{equation}\label{ht-Q}
h_t=h(Q-1).
\end{equation}
Differentiating \eqref{Robin-h} in time yields
\begin{equation}\label{Q-Neumann}
Q_\mu=0
\qquad\hbox{on }\partial\Ct.
\end{equation}

\begin{lemma}[Evolution of the speed ratio]\label{lem-Q}
Let
\[
\cL=\Theta S^{ij}\nabla_i\nabla_j.
\]
Then
\begin{equation}\label{Q-evolution}
(\partial_t-\cL)Q
=
2\Theta S^{ij}(\log h)_iQ_j
+
\bigl(k+1-a_\varphi(\ell/h)\bigr)Q(Q-1).
\end{equation}
\end{lemma}

\begin{proof}
From \eqref{ht-Q},
\[
W_t
=
(Q-1)W+h\nabla^2Q
+\nabla h\otimes\nabla Q+\nabla Q\otimes\nabla h.
\]
Using \eqref{Euler-S},
\[
S_t
=
kS(Q-1)
+hS^{ij}Q_{ij}
+2S^{ij}h_iQ_j.
\]
On the other hand
\[
\partial_t\log\frac{\Theta}{h}
=
\partial_t\log\frac{h\varphi(\ell/h)}f
=
(1-a_\varphi(\ell/h))(Q-1).
\]
Since $Q=\Theta S/h$, combining these identities gives \eqref{Q-evolution}.
\end{proof}

\section{A priori estimates and Orlicz structure}\label{sec-apriori}
\subsection{Uniform $C^0$ and $C^1$ estimates}

Throughout this section, we assume that the solution remains strictly convex on the time interval under
consideration.

\subsubsection{$C^0$ estimate}

Put $u=h/\ell$. By \eqref{Btheta},
\begin{equation*}
u_t
=
\frac{\ell u^2\varphi(1/u)}f\,
\bsig_k(B_\theta[u])-u,
\qquad
u_\mu=0,
\end{equation*}
where $B_\theta[u]$ denotes the right-hand side of \eqref{Btheta}.

Define
\begin{equation*}
g(s)=s^{k+1}\varphi(1/s).
\end{equation*}
Then
\[
\frac{g'(s)}{g(s)}
=
\frac{k+1-a_\varphi(1/s)}s
\le-\frac\delta s.
\]
Hence $g$ is strictly decreasing,
\begin{equation}\label{g-limits}
g(s)\longrightarrow\infty\quad(s\downarrow0),
\qquad
g(s)\longrightarrow0\quad(s\uparrow\infty).
\end{equation}

\begin{proposition}[Uniform support bound]\label{prop-C0}
There are constants $0<c_0<C_0<\infty$, depending only on the fixed data and $u_0$, such that
\[
c_0\le u(\xi,t)\le C_0
\]
on every smooth strictly convex existence interval.
\end{proposition}

\begin{proof}
Let $M(t)=\max_{\overline\Ct}u(\cdot,t)$. At a maximum point, including a boundary maximum,
the condition \eqref{Neumann-u} and the one-sided second-derivative test give
$\nabla u=0$ and $\nabla^2u\le0$. Thus
\[
B_\theta[u]\le M(t)\sig,
\qquad
\bsig_k(B_\theta[u])\le M(t)^k.
\]
Therefore
\[
D^+M
\le
M\left(
\frac{\ell_{\max}}{f_{\min}}g(M)-1
\right).
\]
By \eqref{g-limits}, the bracket is negative when $M$ is sufficiently large.

Similarly, at a minimum point,
\[
B_\theta[u]\ge m(t)\sig,
\qquad
\bsig_k(B_\theta[u])\ge m(t)^k,
\]
and
\[
D_+m
\ge
m\left(
\frac{\ell_{\min}}{f_{\max}}g(m)-1
\right),
\]
which is positive when $m$ is sufficiently small. This gives the claimed time-independent bounds.
\end{proof}

Since $\ell$ has fixed positive upper and lower bounds,
\begin{equation}\label{h-C0}
0<c\le h\le C.
\end{equation}

\subsubsection{$C^1$ estimate}

The strict positivity of $B_\theta[u]$ and \eqref{Btheta} imply, for every unit vector $v$,
\[
\nabla^2u(v,v)
\ge
-C(1+|\nabla u|),
\]
where $C$ depends only on the fixed cap geometry and the $C^0$ bound.

\begin{proposition}[Uniform gradient bound]\label{prop-C1}
There exists $C_1$ such that
\[
|\nabla u|\le C_1
\]
on every smooth strictly convex existence interval.
\end{proposition}

\begin{proof}
Fix a time and suppose $|\nabla u(x_0)|$ is large. On the noncritical set, consider the unit
gradient curve
\[
\gamma'(s)=\frac{\nabla u}{|\nabla u|}.
\]
The Neumann condition makes this vector field tangent to the boundary, so its integral curve remains
in $\overline\Ct$. Put $p(s)=|\nabla u(\gamma(s))|$. Then
\[
p'(s)=\nabla^2u(\gamma',\gamma')
\ge-C(1+p(s)).
\]
For a fixed sufficiently small $L_0>0$, if $|\nabla u(x_0)|$ is large then
$p(s)\ge |\nabla u(x_0)|/4$ on $[0,L_0]$. Hence
\[
\operatorname{osc}_{\Ct}u
\ge
u(\gamma(L_0))-u(\gamma(0))
=
\int_0^{L_0}p(s)\dd s
\ge \frac{L_0}{4}|\nabla u(x_0)|,
\]
contradicting Proposition~\ref{prop-C0}. Thus $|\nabla u|$ is uniformly bounded.
\end{proof}

Consequently
\begin{equation*}
\|h(\cdot,t)\|_{C^1(\overline\Ct)}\le C.
\end{equation*}

\subsection{Uniform control of the speed ratio}

\begin{proposition}[Two-sided $Q$ estimate]\label{prop-Q}
There exist positive constants $c_Q,C_Q$ such that
\[
0<c_Q\le Q(\xi,t)\le C_Q
\]
on every smooth strictly convex existence interval.
\end{proposition}

\begin{proof}
At a spatial maximum of $Q$, the drift in \eqref{Q-evolution} vanishes. If $Q>1$, then
\[
D^+Q_{\max}
\le
-\delta Q_{\max}(Q_{\max}-1).
\]
At a spatial minimum with $Q<1$,
\[
D_+Q_{\min}
\ge
\delta Q_{\min}(1-Q_{\min}).
\]
The boundary condition \eqref{Q-Neumann} permits the same maximum-principle argument at the
boundary. Hence
\[
\min\{1,\min Q(\cdot,0)\}
\le Q(\cdot,t)\le
\max\{1,\max Q(\cdot,0)\}.
\]
\end{proof}

By \eqref{Q-def}, \eqref{h-C0}, and the compact range of $\ell/h$,
\begin{equation}\label{Sk-bounds}
0<c_S\le \bsig_k(W[h])\le C_S.
\end{equation}

\subsection{The dynamic Orlicz square completion}

Let
\[
r=\frac{\ell}{h},
\qquad
\rho=\log\frac h\ell,
\qquad
q_0=\log\ell,
\qquad
H=f\ell^{k-1},
\qquad
K=\log H.
\]
Write
\[
a=a_\varphi(r),\qquad
d=r a_\varphi'(r),\qquad
b=a-2,
\]
and
\[
A_\varphi^{\rm fl}
=
b+d+\frac{b^2}{k+1}.
\]
Under \eqref{supercritical}--\eqref{aphi-monotone},
$A_\varphi^{\rm fl}>0$.

\begin{lemma}[Dynamic Orlicz square completion]\label{lem-square}
Let $e$ be a unit eigenvector of $W[h]$. Then
\begin{equation}\label{square-identity}
\begin{aligned}
&\frac1\Theta
\left[
(k+1)\Theta+\Theta_{ee}
-\frac{k}{k+1}\frac{\Theta_e^2}{\Theta}
+(1-k)\Theta\frac{W(e,e)}h
\right]
\\
={}&
A_\varphi^{\rm fl}(r)
\left[
\rho_e+
\frac{b(r)}
{(k+1)A_\varphi^{\rm fl}(r)}K_e
\right]^2
+
\cT_{\varphi,f}(r)(e,e)
\\
&+
\bigl(3-k-a_\varphi(r)\bigr)\frac{W(e,e)}h,
\end{aligned}
\end{equation}
where $\cT_{\varphi,f}$ is defined in \eqref{Tflow-def}.
\end{lemma}

\begin{proof}
We give the algebra in detail because the coefficient shifts relative to the static problem are
important. Put $L=\log\Theta$. Since
\[
\log h=\rho+q_0,\qquad
\log f=K-(k-1)q_0,
\qquad
r=e^{-\rho},
\]
we have
\begin{equation}\label{L-expression}
L
=
2\rho+(k+1)q_0+\log\varphi(e^{-\rho})-K.
\end{equation}
Differentiating \eqref{L-expression} in the unit direction $e$ gives
\begin{equation}\label{L-first}
L_e=-b\,\rho_e+(k+1)(q_0)_e-K_e.
\end{equation}
Since $r_e=-r\rho_e$, one has
\[
b_e=a_\varphi'(r)r_e=-d\,\rho_e,
\]
and therefore
\begin{equation*}
L_{ee}
=
d\rho_e^2-b\rho_{ee}
+(k+1)(q_0)_{ee}-K_{ee}.
\end{equation*}

We next express $\rho_{ee}$ in geometric variables. From
$W[h](e,e)=h_{ee}+h$ and $W[\ell](e,e)=1$,
\[
(\log h)_{ee}
=
\frac{W(e,e)}h-1-(\rho_e+(q_0)_e)^2,
\]
while
\[
(q_0)_{ee}
=
\frac{\ell_{ee}}\ell-(q_0)_e^2
=
\frac1\ell-1-(q_0)_e^2.
\]
Subtracting yields
\begin{equation}\label{rho-second}
\rho_{ee}
=
\frac{W(e,e)}h-\frac1\ell
-\rho_e^2-2(q_0)_e\rho_e.
\end{equation}

Because $\Theta_{ee}/\Theta=L_{ee}+L_e^2$, substituting
\eqref{L-first}--\eqref{rho-second} into the left-hand side of
\eqref{square-identity} and collecting terms gives
\begin{align}
&(k+1)+L_{ee}+\frac1{k+1}L_e^2
+(1-k)\frac{W(e,e)}h \notag\\
={}&
\left(
b+d+\frac{b^2}{k+1}
\right)\rho_e^2
+\frac{2b}{k+1}\rho_eK_e \notag\\
&-K_{ee}-2(q_0)_eK_e
+\frac1{k+1}K_e^2
+\frac{b+k+1}{\ell}
+\bigl(3-k-a\bigr)\frac{W(e,e)}h.
\label{E-expanded}
\end{align}
The cancellation of the $\rho_e(q_0)_e$ and $(q_0)_e^2$ terms is the reason for introducing
$H=f\ell^{k-1}$.

Completing the square in $\rho_e$ gives
\begin{align}
&A_\varphi^{\rm fl}\rho_e^2
+\frac{2b}{k+1}\rho_eK_e+\frac1{k+1}K_e^2 \notag\\
&\qquad=
A_\varphi^{\rm fl}\left(
\rho_e+\frac{b}{(k+1)A_\varphi^{\rm fl}}K_e
\right)^2
+
\frac{b+d}{(k+1)A_\varphi^{\rm fl}}K_e^2.
\label{rho-square}
\end{align}
Combining \eqref{E-expanded} and \eqref{rho-square} gives exactly
\eqref{square-identity} and the tensor formula \eqref{Tflow-def}.

For later use we record an equivalent factorization. For a fixed parameter $s>0$, set
\begin{equation*}
c_\varphi^{\rm fl}(s)=\frac{b(s)^2}{b(s)+d(s)},
\qquad
\alpha_\varphi^{\rm fl}(s)
=
\frac1{k+1+c_\varphi^{\rm fl}(s)}.
\end{equation*}
With $s$ held fixed, $\alpha_\varphi^{\rm fl}(s)$ is constant under spatial differentiation, and
\[
\frac{
\nabla^2(\ell H^{-\alpha_\varphi^{\rm fl}(s)})
+\ell H^{-\alpha_\varphi^{\rm fl}(s)}\sig
}{
\ell H^{-\alpha_\varphi^{\rm fl}(s)}
}
=
\frac1\ell\sig
-\alpha_\varphi^{\rm fl}(s)\nabla^2K
-2\alpha_\varphi^{\rm fl}(s)dq_0\od dK
+\bigl(\alpha_\varphi^{\rm fl}(s)\bigr)^2dK\otimes dK.
\]
Using
\[
\alpha_\varphi^{\rm fl}(s)
=
\frac{b(s)+d(s)}
{(k+1)(b(s)+d(s))+b(s)^2}
=
\frac{b(s)+d(s)}
{(k+1)A_\varphi^{\rm fl}(s)}
\]
and
\[
b(s)-c_\varphi^{\rm fl}(s)
=
\frac{b(s)d(s)}{b(s)+d(s)},
\]
we obtain
\begin{equation*}
\begin{aligned}
\cT_{\varphi,f}(s)
={}&
\bigl(k+1+c_\varphi^{\rm fl}(s)\bigr)
\frac{
\nabla^2(\ell H^{-\alpha_\varphi^{\rm fl}(s)})
+\ell H^{-\alpha_\varphi^{\rm fl}(s)}\sig
}{
\ell H^{-\alpha_\varphi^{\rm fl}(s)}
}\\
&+
\frac{b(s)d(s)}{\ell\,[b(s)+d(s)]}\,\sig .
\end{aligned}
\end{equation*}
After evaluating this identity at $s=r=\ell/h$, we see that the condition $d\ge0$ and a strict
spherical-convexity lower bound for the first term imply \eqref{dynamic-convexity}.
\end{proof}

\begin{remark}[Power-law check]\label{rem-power-check}
If $\varphi(s)=s^p$, then $a_\varphi=p$, $d=0$, and
\[
A_\varphi^{\rm fl}
=
\frac{(p-2)(p+k-1)}{k+1}.
\]
A direct computation gives
\[
\cT_{\varphi,f}
=
(p+k-1)
\frac{
\nabla^2(\ell H^{-1/(p+k-1)})
+\ell H^{-1/(p+k-1)}\sig
}{
\ell H^{-1/(p+k-1)}
}.
\]
Thus the dynamic structural tensor recovers exactly the exponent
$1/(p+k-1)$ appearing in the closed flow of Ivaki \cite{Ivaki2019}.
\end{remark}

\section{Convexity and curvature estimates}\label{sec-curvature}
\subsection{Preservation of strict convexity}

In this section the strict inequality \eqref{dynamic-convexity} is not needed; it suffices that
\begin{equation*}
\cT_{\varphi,f}(s)\ge0\qquad(s>0).
\end{equation*}

Define
\[
q=h\,\lambda_{\min}(W[h]).
\]
The smallest eigenvalue is only Lipschitz in general, so the correct formulation is in the viscosity
sense, following \cite{BISConstantRank,Ivaki2019}.

\begin{lemma}[Viscosity inequality for the minimum radius]\label{lem-q-viscosity}
% R04 | Positive parabolic lower-support domain
% Status: verified clarification
Let $h$ be a smooth solution with $h>0$ and $W[h]>0$ near an interior
space-time point $(x_0,t_0)$, where $t_0>0$.  Let $\xi>0$ be a smooth
parabolic lower support of $q=h\lambda_{\min}(W)$ there: $\xi\le q$
in a neighborhood of $x_0$ for times immediately preceding $t_0$, with
equality at $(x_0,t_0)$.  If $e$ is a unit minimum eigenvector at the
contact point, then

\begin{equation}\label{q-viscosity}
\begin{aligned}
(\partial_t-\cL)\xi
&+
2\Theta S^{ij}\xi_i
\nabla_j\log\frac{\xi}{h}
\\
&+
\bigl[(a_\varphi(r)+k-3)Q+2\bigr]\xi
\ge
h^2Q\,\cT_{\varphi,f}(r)(e,e).
\end{aligned}
\end{equation}
\end{lemma}

% R05 | Complete minimum-eigenvalue proof prefix
% Status: independent analytic derivations agree and explicit remainder checked; no extension to nonpositive W
\begin{proof}
All quantities below are evaluated at the interior contact point
$(x_0,t_0)$.  Choose an orthonormal frame which is geodesic at $x_0$ and
diagonalizes $W$, and write
\[
W_{ij}=\lambda_i\delta_{ij},\qquad
\lambda:=\lambda_1=\cdots=\lambda_m<\lambda_{m+1}\le\cdots\le\lambda_n,
\qquad e=e_1.
\]
Put
\[
A^{ij}=\Theta S^{ij},\qquad
D_i=\frac{\xi_i-\lambda h_i}{h},\qquad
\eta_{ij}=\nabla_1W_{ij}.
\]
For each unit vector in the minimum eigenspace, $hW(e,e)-\xi$ has a
spatial minimum of zero.  Differentiating and polarizing within that
eigenspace gives
\begin{equation}\label{minimum-contact-first}
\nabla_iW_{ab}=\delta_{ab}D_i\qquad(1\le a,b\le m).
\end{equation}
The second-derivative test, allowing the test eigenvector to vary in the
complementary eigenspaces, gives
\begin{equation}\label{minimum-contact-second}
\xi_{ii}\le (hW_{11})_{;ii}
 -2h\sum_{p>m}\frac{(\nabla_iW_{1p})^2}{\lambda_p-\lambda}.
\end{equation}
Indeed, one obtains this formula by choosing the $p$-component of the
first derivative of the unit test vector to be
$-\nabla_iW_{1p}/(\lambda_p-\lambda)$ and completing the resulting
quadratic expression.  No division by a gap inside the minimum
eigenspace is needed.  The parabolic lower-support condition gives
$\xi_t\ge\partial_t(hW_{11})$: the difference is nonnegative for past
times and vanishes at the terminal contact time.

Since $W$ is Codazzi,
\[
\eta_{ip}=\nabla_1W_{ip}=\nabla_iW_{1p},\qquad
\eta_{ip}=\delta_{1p}D_i\quad(p\le m).
\]
At a diagonal matrix, write $S_i=S^{ii}>0$.  The flow, Euler homogeneity,
and the unit-sphere commutator yield
\begin{align*}
(\partial_t-\cL)h
 &= (1-k)\Theta S-h+\Theta h\sum_iS_i,\\
(\partial_t-\cL)W_{11}
 &=S\Theta_{11}+2\Theta_1S_1^{\mathrm{dir}}
   +\Theta S^{ij,rs}\eta_{ij}\eta_{rs}
   +(k+1)\Theta S
   -\left(\Theta\sum_iS_i+1\right)\lambda,
\end{align*}
where $S_1^{\mathrm{dir}}:=\nabla_1S=S^{ij}\eta_{ij}$; this notation
distinguishes the directional derivative from $S_1=S^{11}$.
Using \eqref{minimum-contact-second} in the product rule, and using
$\Theta S=hQ$ and $\xi=h\lambda$, gives
\begin{align}
&(\partial_t-\cL)\xi
 +2\Theta S^{ij}\xi_i\nabla_j\log\frac{\xi}{h}
 +\bigl[(k-1)Q+2\bigr]\xi \notag\\
&\quad\ge
 hS\bigl[(k+1)\Theta+\Theta_{11}\bigr]
 +2h\Theta_1S_1^{\mathrm{dir}}
 +h\Theta\Bigg[
 S^{ij,rs}\eta_{ij}\eta_{rs}
 +2\sum_iS_i\left(
     \frac{D_i^2}{\lambda}
     +\sum_{p>m}\frac{\eta_{ip}^2}{\lambda_p-\lambda}
                  \right)\Bigg].
\label{minimum-product-expansion}
\end{align}
In particular, the drift is exact: since
$\nabla_i\log(\xi/h)=D_i/\lambda$, adding it cancels the product term
$-2\Theta S^{ij}h_iD_j$ and leaves
$2h\Theta\sum_iS_iD_i^2/\lambda$.

We record explicitly the inverse-concavity inequality used here.  For
$W>0$ and every symmetric matrix $\eta$,
\begin{equation}\label{inverse-concavity-S}
\left(S^{ij,rs}+2S^{ir}(W^{-1})^{js}\right)\eta_{ij}\eta_{rs}
\ge\frac{k+1}{kS}(S^{ij}\eta_{ij})^2.
\end{equation}
Here is a direct verification, including its domain.  Set $F=S^{1/k}$
and $F_*(B)=1/F(B^{-1})$ on the positive definite cone.  On positive
coordinate vectors $x$, let $\mathcal I$ be the collection of $k$-element
subsets of $\{1,\ldots,n\}$ and put
\[
g_I(x)=\left(\prod_{i\in I}x_i\right)^{1/k},\qquad
P(y)=\left(\sum_{I\in\mathcal I}y_I^{-k}\right)^{-1/k}.
\]
Each $g_I$ is positive and concave.  The function $P$ is increasing in
each variable and concave: writing $Z=\sum_Iy_I^{-k}$, one has
\[
D^2P(y)[v,v]
=(k+1)Z^{-1/k-2}
 \left[
  \left(\sum_Iy_I^{-k-1}v_I\right)^2
  -Z\sum_Iy_I^{-k-2}v_I^2
 \right]\le0
\]
by Cauchy--Schwarz.  Consequently
\[
F_*(x)=\binom nk^{1/k}P\bigl((g_I(x))_{I\in\mathcal I}\bigr)
\]
is a symmetric concave function of the positive eigenvalues.  It is
also concave as a spectral function on positive definite matrices.
For completeness, at distinct eigenvalues its second derivative is
\[
D^2F_*(B)[Z,Z]
 =\sum_{i,j}(F_*)_{ij}Z_{ii}Z_{jj}
  +2\sum_{i<j}\frac{(F_*)_i-(F_*)_j}{x_i-x_j}Z_{ij}^2\le0;
\]
the first term is nonpositive by scalar concavity, and the quotients
are nonpositive by symmetry and concavity.  The case of repeated eigenvalues follows
by continuity.

Now differentiate $F_*$ along
$B(s)=W^{-1}-sW^{-1}\eta W^{-1}$, which remains positive definite for
small $|s|$.  If $A(s)=B(s)^{-1}$, then
$A'(0)=\eta$ and $A''(0)=2\eta W^{-1}\eta$.  The inequality
$\frac{d^2}{ds^2}F_*(B(s))|_{s=0}\le0$ is therefore
\[
F^{ij,rs}\eta_{ij}\eta_{rs}
 +2F^{ir}(W^{-1})^{js}\eta_{ij}\eta_{rs}
\ge\frac2F(F^{ij}\eta_{ij})^2.
\]
Substituting $F=S^{1/k}$ and multiplying by $kS^{1-1/k}$ gives
\eqref{inverse-concavity-S}, with the coefficient
$2/k-(1/k-1)=(k+1)/k$.

The bracket in \eqref{minimum-product-expansion} is bounded below by
the left-hand side of \eqref{inverse-concavity-S}.  Indeed, its difference
from that expression is exactly
\[
2\sum_iS_i\sum_{p>m}
 \left(\frac1{\lambda_p-\lambda}-\frac1{\lambda_p}\right)
 \eta_{ip}^2\ge0,
\]
because \eqref{minimum-contact-first} supplies all the terms with
$p\le m$ and $0<\lambda<\lambda_p$.  Thus
\[
h\Theta\frac{k+1}{kS}(S_1^{\mathrm{dir}})^2
 +2h\Theta_1S_1^{\mathrm{dir}}
\ge-hS\frac{k}{k+1}\frac{\Theta_1^2}{\Theta}.
\]
Inserting this into \eqref{minimum-product-expansion} proves
\begin{align}
(\partial_t-\cL)\xi
&+2\Theta S^{ij}\xi_i\nabla_j\log\frac{\xi}{h}
 +\bigl[(k-1)Q+2\bigr]\xi \notag\\
&\ge hS\left[(k+1)\Theta+\Theta_{ee}
                  -\frac{k}{k+1}\frac{\Theta_e^2}{\Theta}\right].
\label{visc-intermediate}
\end{align}

We now insert the Orlicz coefficient. Lemma~\ref{lem-square} implies
\begin{align}
&\frac1\Theta\left[
(k+1)\Theta+\Theta_{ee}
-\frac{k}{k+1}\frac{\Theta_e^2}{\Theta}
\right]\notag\\
&\qquad=
A_\varphi^{\rm fl}
\left[
\rho_e+\frac{b}{(k+1)A_\varphi^{\rm fl}}K_e
\right]^2
+\cT_{\varphi,f}(e,e)
+(2-a_\varphi(r))\frac{W(e,e)}h.
\label{visc-square-insert}
\end{align}
Indeed, \eqref{visc-square-insert} is obtained from \eqref{square-identity} by subtracting
$(1-k)W(e,e)/h$ from both sides. At the contact point
\[
hW(e,e)=\xi,
\qquad
\Theta S=hQ.
\]
Multiplying \eqref{visc-square-insert} by $hS\Theta=h^2Q$, discarding the nonnegative square,
and moving
\[
(2-a_\varphi(r))Q\,\xi
\]
to the left of \eqref{visc-intermediate} yields
\[
(\partial_t-\cL)\xi
+
2\Theta S^{ij}\xi_i\nabla_j\log\frac{\xi}{h}
+
\bigl[(a_\varphi(r)+k-3)Q+2\bigr]\xi
\ge
h^2Q\,\cT_{\varphi,f}(r)(e,e),
\]
which is \eqref{q-viscosity}.
\end{proof}

\begin{lemma}[Boundary derivative of a minimum-radius support]\label{lem-q-boundary}
Assume
\[
0<\theta\le\frac{\pi}{2},
\qquad
W>0,
\qquad
H_\mu\ge0.
\]
At a boundary point, choose a unit minimum eigendirection adapted to the boundary splitting:
if the minimum eigenspace contains a tangential direction, choose such a direction; otherwise
choose $e=\mu$.  Extend $e$ smoothly near the point and set $v=hW(e,e)$.  Then $v$ is a
smooth local upper support of $q=h\lambda_{\min}(W)$ and
\[
v_\mu+(n-1)\cot\theta\,v\ge0.
\]
\end{lemma}

\begin{proof}
At the boundary choose an orthonormal frame
$e_1,\ldots,e_{n-1},e_n=\mu$ which diagonalizes the tangential block of $W$.
The mixed components vanish by \eqref{W-mixed}.

If the minimum eigenvalue is tangential, say $\lambda_\alpha$, then
\[
v_\mu
=
h_\mu\lambda_\alpha+hW_{\alpha\alpha;\mu}
=
h\cot\theta\,\lambda_n\ge0.
\]
For $0<\theta<\pi/2$ the sign is strict, while at $\theta=\pi/2$ it is zero.

Suppose now that the minimum eigenvalue is the normal one $\lambda_n$. If $k=1$, then
$Q_\mu=0$ and $(\ell/h)_\mu=0$ give
\[
\frac{(\sigma_1(W))_\mu}{\sigma_1(W)}
=
\frac{f_\mu}{f}-\cot\theta.
\]
Using
\[
\sigma_1(W)=\sum_{\alpha<n}\lambda_\alpha+\lambda_n,
\qquad
(\lambda_\alpha)_\mu
=
\cot\theta(\lambda_n-\lambda_\alpha),
\]
one obtains directly
\[
\frac{v_\mu}{v}
=
\left(
1+\frac{\sum_{\alpha<n}\lambda_\alpha}{\lambda_n}
\right)\frac{f_\mu}{f}
-(n-1)\cot\theta.
\]
Here $H=f$, so $H_\mu\ge0$ gives
\[
\frac{v_\mu}{v}\ge-(n-1)\cot\theta.
\]

Assume next that $2\le k\le n-1$ and write
\[
s_j=\sigma_j(\lambda_1,\ldots,\lambda_{n-1}).
\]
Since $Q_\mu=0$ and $(\ell/h)_\mu=0$,
\[
\frac{(\sigma_k(W))_\mu}{\sigma_k(W)}
=
\frac{f_\mu}{f}-\cot\theta.
\]
Using
\[
\sigma_k(W)=s_k+\lambda_ns_{k-1},
\]
the boundary identity \eqref{W-normal-derivative}, and the identities
\[
\sum_{\alpha<n}\sigma_{k-1}(\lambda'|\alpha)=(n-k)s_{k-1},
\]
\[
\sum_{\alpha<n}\lambda_\alpha\sigma_{k-1}(\lambda'|\alpha)=ks_k,
\]
\[
\sum_{\alpha<n}\sigma_{k-2}(\lambda'|\alpha)=(n-k+1)s_{k-2},
\]
\[
\sum_{\alpha<n}\lambda_\alpha\sigma_{k-2}(\lambda'|\alpha)=(k-1)s_{k-1},
\]
% R06 | Normal boundary algebra for k>=2
% Status: verified local derivation; k=1 remains separately treated
we may calculate as follows. Put $c_\theta=\cot\theta$ and $r_n=\lambda_n$.
The tangential derivative formula gives
\[
(s_j)_\mu
=c_\theta\bigl[r_n(n-j)s_{j-1}-js_j\bigr]
\qquad(1\le j\le n-1),
\]
where $s_0=1$.  Hence
\begin{align*}
(\sigma_k(W))_\mu
={}&(r_n)_\mu s_{k-1}\\
&+c_\theta\bigl[
 (n-2k+1)r_ns_{k-1}-ks_k
 +(n-k+1)r_n^2s_{k-2}\bigr].
\end{align*}
Equating this with
$(s_k+r_ns_{k-1})(f_\mu/f-c_\theta)$ and using
$v_\mu/v=c_\theta+(r_n)_\mu/r_n$ gives
\[
\frac{v_\mu}{v}
=\left(1+\frac{s_k}{r_ns_{k-1}}\right)\frac{f_\mu}{f}
+c_\theta\left[
 (k-1)\frac{s_k}{r_ns_{k-1}}-(n-2k+1)
 -(n-k+1)r_n\frac{s_{k-2}}{s_{k-1}}\right].
\]
Finally, substitute $f_\mu/f=H_\mu/H-(k-1)c_\theta$ to obtain
\begin{equation}\label{q-normal-boundary}
\frac{v_\mu}{v}
=
\left(
1+\frac{s_k}{\lambda_ns_{k-1}}
\right)\frac{H_\mu}{H}
-
\left(
n-k+(n-k+1)\lambda_n\frac{s_{k-2}}{s_{k-1}}
\right)\cot\theta.
\end{equation}
Moreover,
\[
(k-1)s_{k-1}
=
\sum_{\alpha<n}
\lambda_\alpha\sigma_{k-2}(\lambda'|\alpha)
\ge
\lambda_n(n-k+1)s_{k-2},
\]
so
\[
(n-k+1)\lambda_n\frac{s_{k-2}}{s_{k-1}}\le k-1.
\]
Hence $H_\mu\ge0$ and $\cot\theta\ge0$ imply
\[
\frac{v_\mu}{v}\ge-(n-1)\cot\theta.
\]
This proves the claim in all cases.
\end{proof}

% R07 | Finite-time convexity barrier, both boundary directions
% Status: verified local completion; incorporates minor text correction
\begin{proof}[Proof of Theorem~\ref{thm-preserve}]
Fix $T>0$ up to which the solution is smooth.  Choose a smooth function
$\beta$ on $\overline\Ct$ such that $\beta\ge0$, $\beta=0$ on the
boundary, and $\beta_\mu=-1$ there.  Fix $\zeta>0$ and set
\[
\psi=-(n-1)\log\ell+\zeta\beta,
\qquad
\xi(x,t)=\varepsilon e^{-Mt+\psi(x)}.
\]
Choose
\[
0<\varepsilon<\min_{\overline\Ct}
                 \bigl(q(x,0)e^{-\psi(x)}\bigr).
\]
Thus $q(x,0)>\xi(x,0)$ for every $x\in\overline\Ct$.

We work up to a possible first loss of positive definiteness of $W$.
On this interval, the support bounds and smoothness on $[0,T]$ give
upper bounds for $A^{ij}=\Theta S^{ij}$, $\nabla h/h$, and
\[
\mathfrak b=(a_\varphi(\ell/h)+k-3)Q+2.
\]
These bounds require no positive lower bound for the eigenvalues of
$W$: the derivatives $S^{ij}$ are polynomials in its entries.  Since
$\psi$ is fixed, there is a finite constant $C_T$, independent of
$M$ and $\varepsilon$, such that
\[
\left|
 A^{ij}(-\psi_{ij}+\psi_i\psi_j-2\psi_i h_j/h)
 +\mathfrak b
\right|\le C_T
\]
throughout the positive-definite interval.  Fix $M>C_T$.

Suppose that $q-\xi$ has a first zero at $(x_*,t_*)$ with $t_*>0$.
At this contact point, $q=\xi>0$, so $W$ is still positive definite.
If the contact is interior, $\xi$ is a positive parabolic lower
support and Lemma~\ref{lem-q-viscosity} applies.  Direct differentiation
of $\xi$ gives
\begin{align*}
&(\partial_t-\cL)\xi
 +2\Theta S^{ij}\xi_i\nabla_j\log\frac{\xi}{h}
 +\mathfrak b\xi\\
&\quad=
\xi\left[-M+
 A^{ij}(-\psi_{ij}+\psi_i\psi_j-2\psi_i h_j/h)
 +\mathfrak b\right]<0.
\end{align*}
This contradicts the nonnegative right-hand side of
\eqref{q-viscosity}, since $\cT_{\varphi,f}\ge0$.

Suppose instead that the first contact lies on the boundary.
The identity $W_{\alpha n}=0$ splits $W$ into its tangential block and
its normal entry.  If the minimum eigenspace contains a tangential
unit vector, choose one; otherwise, that eigenspace is the normal line,
and we choose $e=\mu$.  Extend $e$ to a smooth unit vector field near the
point and put $v=hW(e,e)$.  Then $v\ge q\ge\xi$, with equality at the
contact point, and therefore $(v-\xi)_\mu\le0$ for the outward conormal.
In both cases Lemma~\ref{lem-q-boundary} gives
\[
v_\mu+(n-1)\cot\theta\,v\ge0.
\]
On the other hand, $\ell_\mu/\ell=\cot\theta$ and
$\beta_\mu=-1$ imply
\[
\xi_\mu+(n-1)\cot\theta\,\xi=-\zeta\xi<0.
\]
Since $v=\xi$ at the contact point, these inequalities contradict
$(v-\xi)_\mu\le0$.

Thus no first contact is possible.  In particular,
$q\ge\varepsilon e^{-MT+\min\psi}>0$ up to any proposed first loss
of positive definiteness.  Continuity excludes such a loss in
$[0,T]$.  Hence $W[h]>0$ throughout $[0,T]$.
\end{proof}

\subsection{A time-uniform non-collapsing estimate}

We now assume the strict dynamic condition \eqref{dynamic-convexity} and the quantitative boundary
condition \eqref{boundary-quant}.

\begin{lemma}[Quantitative normal boundary sign]\label{lem-normal-quant}
At a boundary point at which the normal radius $\lambda_n$ is the minimum eigenvalue, the
normal support $v=hW(\mu,\mu)$ satisfies
\[
v_\mu\ge0.
\]
\end{lemma}

\begin{proof}
If $k=1$, then $\lambda_\alpha\ge\lambda_n$ gives
\[
1+\frac{\sum_{\alpha<n}\lambda_\alpha}{\lambda_n}\ge n.
\]
The boundary assumption \eqref{boundary-quant} becomes
\[
\frac{f_\mu}{f}\ge\frac{n-1}{n}\cot\theta,
\]
so the $k=1$ formula in Lemma~\ref{lem-q-boundary} gives $v_\mu\ge0$.

Assume $2\le k\le n-1$ and use $s_j$ from Lemma~\ref{lem-q-boundary}. The identity
\[
ks_k
=
\sum_{\alpha<n}
\lambda_\alpha\sigma_{k-1}(\lambda'|\alpha)
\]
and $\lambda_\alpha\ge\lambda_n$ give
\[
\frac{s_k}{\lambda_ns_{k-1}}\ge\frac{n-k}{k}.
\]
Together with
\[
(n-k+1)\lambda_n\frac{s_{k-2}}{s_{k-1}}\le k-1,
\]
this gives the two coefficient bounds needed in \eqref{q-normal-boundary}.
Moreover,
\[
\frac{H_\mu}{H}
=
\frac{f_\mu}{f}+(k-1)\cot\theta
\ge
\frac{k(n-1)}n\cot\theta
\]
by \eqref{boundary-quant}. Substitution into
\eqref{q-normal-boundary} yields
\[
\frac{v_\mu}{v}\ge0.
\]
At $\theta=\pi/2$ all $\cot\theta$ terms vanish, so the same conclusion holds.
\end{proof}

\begin{proposition}[Uniform lower principal-radius bound]\label{prop-lower-radius}
Under the geometric, data, and initial hypotheses of
Theorem~\ref{thm-main}, there exists $c>0$ such that
\[
h\,\lambda_{\min}(W[h])\ge c
\qquad\hbox{for all }t\ge0
\]
as long as the smooth solution exists. This estimate does not require
Lemmas~\ref{lem-parabolic} and \ref{lem-EKH}.
\end{proposition}

% R08 | Uniform non-collapsing without circular constants
% Status: verified local completion
\begin{proof}
The previously established estimates provide constants, independent of
time on the smooth strictly convex existence interval, such that
\[
0<h_-\le h\le h_+,\qquad
0<Q_-\le Q\le Q_+,\qquad
S\le S_+,\qquad
\Theta\le\Theta_+,\qquad
|\nabla h|\le H_1.
\]
Put
\[
c_*:=\tau_0h_-^2Q_->0,
\qquad
B_*:=\sup\left|(a_\varphi(\ell/h)+k-3)Q+2\right|<\infty.
\]
The compact range of $\ell/h$ makes $B_*$ finite.  Thus the right-hand
side of \eqref{q-viscosity} is at least $c_*$.

Let $q_{0,-}=\min_{\overline\Ct}q(\cdot,0)>0$ and choose
\[
0<\varepsilon<
\min\left\{\frac{q_{0,-}}4,\frac{c_*}{8(B_*+1)}\right\}.
\]
Choose a fixed smooth function $\beta\ge0$ with $\beta=0$ and
$\beta_\mu=-1$ on the boundary.  We will use
$\xi(x)=\varepsilon e^{\zeta\beta(x)}$, where $\zeta>0$ will be
chosen after estimating the coefficients in terms of $\varepsilon$.

Before and at any first contact of $q$ with $\xi$, the inequality
$q\ge\xi\ge\varepsilon$ implies
\[
\lambda_{\min}(W)\ge\frac{\varepsilon}{h_+}.
\]
A monomial of $\sigma_k(W)$ containing $\lambda_{\max}$ and any
$k-1$ other eigenvalues therefore yields
\[
\lambda_{\max}(W)\le
L_\varepsilon:=\binom nk S_+
                  \left(\frac{h_+}{\varepsilon}\right)^{k-1}.
\]
This also covers $k=1$, with the empty product equal to one.
Consequently the positive matrix $A^{ij}=\Theta S^{ij}$ satisfies
\[
\|A\|_{\mathrm{op}}\le
\Lambda_\varepsilon:=\Theta_+\frac{k}{n}L_\varepsilon^{k-1}.
\]
This follows from
$S^{ii}=\binom nk^{-1}\sigma_{k-1}(\lambda|i)$ at a diagonal matrix.
All these bounds are independent of the proposed contact time and of
$\zeta$.

Write $K_1=\|\nabla\beta\|_\infty$,
$K_2=\|\nabla^2\beta\|_{\mathrm{op},\infty}$ and
$K_0=\|\beta\|_\infty$, and set
\[
D_\varepsilon
=n\Lambda_\varepsilon
 \left(K_2+2K_1H_1/h_-+K_1^2\right).
\]
Choose
\[
0<\zeta\le
\min\left\{1,\frac{\log2}{K_0+1},
                 \frac{c_*}{8\varepsilon(D_\varepsilon+1)}\right\}.
\]
Then $\varepsilon\le\xi\le2\varepsilon<q_{0,-}$, so the initial
comparison holds at every point of $\overline\Ct$.
At a proposed contact, direct differentiation gives
\begin{align*}
&-\cL\xi
 +2\Theta S^{ij}\xi_i\nabla_j\log\frac{\xi}{h}\\
&\qquad=
\xi\zeta A^{ij}
 \left(-\beta_{ij}+\zeta\beta_i\beta_j-2\beta_i h_j/h\right).
\end{align*}
Its absolute value is at most
$\xi\zeta D_\varepsilon\le c_*/4$.
The absolute value of the zeroth-order contribution is at most
$B_*\xi\le2B_*\varepsilon<c_*/4$.
Since $\xi_t=0$, an interior first contact would therefore give
$c_*\le c_*/2$ in \eqref{q-viscosity}, a contradiction.

At a first contact on the boundary, use the splitting $W_{\alpha n}=0$ to
choose a tangential minimum eigenvector if one exists and $e=\mu$
otherwise.  Extend $e$ as a smooth unit vector field and set
$v=hW(e,e)$.  Then $v\ge q\ge\xi$ and $v=\xi$ at contact, so
$(v-\xi)_\mu\le0$.  For a tangential minimum, the boundary identity
gives $v_\mu=h\cot\theta\,\lambda_n\ge0$.  For a normal minimum,
Lemma~\ref{lem-normal-quant} gives the same inequality.  But
\[
\xi_\mu=-\zeta\xi<0,
\]
which contradicts $(v-\xi)_\mu\le0$ in either case, including
$\theta=\pi/2$.

There is therefore no first contact.  Since all constants were chosen
from the initial data and the preceding time-uniform estimates,
\[
h\lambda_{\min}(W)=q\ge\xi\ge\varepsilon>0
\]
on the entire smooth existence interval.
\end{proof}

Combining Proposition~\ref{prop-lower-radius} with \eqref{h-C0} gives
\begin{equation}\label{lambda-lower}
\lambda_{\min}(W)\ge c>0.
\end{equation}

\begin{proposition}[Two-sided radius estimate]\label{prop-C2}
There exist $0<c<C<\infty$ such that
\[
c\sig\le W[h(\cdot,t)]\le C\sig
\]
on the whole existence interval.
\end{proposition}

\begin{proof}
The lower bound is \eqref{lambda-lower}. If $\lambda_n=\lambda_{\max}(W)$, then one term in
$\sigma_k(W)$ contains $\lambda_n$ times $k-1$ other eigenvalues, each bounded below by $c$.
Hence
\[
\sigma_k(W)\ge c^{k-1}\lambda_n.
\]
The upper bound in \eqref{Sk-bounds} therefore yields $\lambda_n\le C$.
\end{proof}

\section{Long-time behavior and contact-angle issues}\label{sec-longtime}
\subsection{Long-time existence}\label{subsec-longtime}

\begin{lemma}[Quantitative interior concave estimate]\label{lem-EKH}
Let $n\geq2$ and let
$\mathfrak F:B_1\times\operatorname{Sym}(n)\to\mathbb R$
be continuous, concave in its matrix variable, and uniformly elliptic:
\[
 \lambda\operatorname{tr}N\leq
 \mathfrak F(x,M+N)-\mathfrak F(x,M)
 \leq\Lambda\operatorname{tr}N\qquad(N\geq0).
\]
Suppose
\[
 |\mathfrak F(x,M)-\mathfrak F(y,M)|
 \leq L(1+\|M\|)|x-y|,
 \qquad |\mathfrak F(x,0)|\leq L.
\]
If $g\in C^\alpha(B_1)$, $0<\alpha<1$, and a bounded classical
solution satisfies $\mathfrak F(x,D^2v)=g$, then
\[
 \|v\|_{C^{2,\beta}(B_{1/2})}
 \leq C\bigl(\|v\|_{L^\infty(B_1)}+1+\|g\|_{C^\alpha(B_1)}\bigr)
\]
for some $0<\beta<\alpha$. The constants depend only on
$n,\lambda,\Lambda,L,\alpha$ and the chosen admissible $\beta$.
No differentiability in the matrix variable is required.
\end{lemma}
\begin{proof}
Set $U=-v$ and
$\mathcal G(x,z,p,M)=-\mathfrak F(x,-M)+g(x)$.
Then $\mathcal G$ is convex in $M$, independent of $(z,p)$,
and satisfies the rank-one ellipticity condition with
$\nu=\min\{1,\lambda,\Lambda^{-1}\}$.
Choose $0<\beta<\min\{\alpha,\bar\alpha(n,\nu)\}$.
On a ball of diameter at most two,
\[
 [\mathcal G(\cdot,z,p,M)]_{C^\beta}
 \leq 2^{1-\beta}L(1+\|M\|)
       +2^{\alpha-\beta}[g]_{C^\alpha}.
\]
Together with $|\mathcal G(x,0,0,0)|\leq L+\|g\|_\infty$,
this verifies all of Safonov's conditions (F0)--(F4)
\cite[Assumptions 3.1, p.~9]{SafonovBVP}, including the required
linear growth of the spatial H\"older modulus in $\sum_{i,j}|M_{ij}|$.
Apply his Theorem~3.1 and estimate (3.6), p.~10, on a ball
compactly contained in $B_1$, with boundary values taken from $U$.
The resulting interior $C^{2,\beta}$ solution, continuous up to
the boundary of the ball, equals $U$ by classical comparison
\cite[Corollary 3.1, p.~12]{SafonovBVP}.
The claim follows by restricting the weighted estimate to a smaller ball
and using a finite covering of $B_{1/2}$.
\end{proof}

\begin{remark}[Choice of the interior estimate]\label{rem-EKH-source}
Theorem~6.2 of Lieberman--Trudinger \cite{LiebermanTrudinger},
under its condition (6.1), requires a spatial derivative bound for
the operator. Applied with right-hand side
$g=(fQ/(h\varphi(\ell/h)))^{1/k}$, that bound would require an
unavailable uniform bound for $Dg$. A spatial $C^\alpha$ bound for
$g$ alone does not verify that application. We instead use
Lemma~\ref{lem-EKH}, which requires only the stated spatial H\"older
norm of the right-hand side.
\end{remark}

\begin{lemma}[Boundary second-derivative H\"older estimate]\label{lem-boundary-holder}
Fix $n\geq2$, $0<\theta<\pi$, $1\leq k\leq n$, $0<c<C$, $M_0<\infty$, and
$0<\alpha<1$. Let $h$ be smooth on the closed spherical cap
$\overline{\mathcal C_\theta}$ and suppose
\[
 h_\mu=\cot\theta\,h,\qquad
 c\sigma\leq W[h]:=\nabla^2h+h\sigma\leq C\sigma,
 \qquad \|h\|_{C^2}\leq M_0.
\]
If
\[
 \bar\sigma_k(W[h])^{1/k}=g,\qquad
 \|g\|_{C^\alpha(\overline{\mathcal C_\theta})}\leq M_1,
\]
then, for some $\beta>0$ depending only on the indicated
fixed structural data,
\[
 \|h\|_{C^{2,\beta}(\overline{\mathcal C_\theta})}\leq C_*.
\]
Here $C_*$ depends on $n,k,\theta,c,C,M_0,M_1,\alpha$ and the constants
in Lemma~\ref{lem-EKH}, and not on higher derivatives of $h$ or $g$.
\end{lemma}

\begin{proof}
Set $u=h/\ell$. The fixed positive function $\ell$ satisfies
$\nabla^2\ell+\ell\sigma=\sigma$ and
$\ell_\mu=\cot\theta\,\ell$, so
\[
 u_\mu=0,\qquad
 W[h]=\ell\nabla^2u+d\ell\otimes du+du\otimes d\ell+u\sigma.
\]
In particular, $\|u\|_{C^2}$ is bounded by fixed data. The second
identity, including the term $u\sigma$, is used below.

Choose a collar of the boundary, let $r\geq0$ be the inward geodesic distance,
and use angular coordinates $y$ on $S^{n-1}$. The round metric and
$\ell$ are
\[
 \sigma=dr^2+\psi(r)^2\gamma,
 \qquad \psi(r)=\sin(\theta-r),
 \qquad \ell(r)=1-\cos\theta\cos(\theta-r),
\]
where $\gamma$ is the unit round metric on $S^{n-1}$. We choose and fix the collar
small enough that $\psi>0$. Let
$E_a=E_a^\alpha(y)\partial_{y^\alpha}$, $a=1,\ldots,n-1$, be a smooth
local $\gamma$-orthonormal frame. An orthonormal frame for $\sigma$ is
$e_a=\psi^{-1}E_a$ and $e_n=\partial_r$. In this frame direct computation
gives
\begin{align*}
 W_{ab}
 &=\frac{\ell}{\psi^2}E_a^\alpha E_b^\beta
   \bigl(u_{\alpha\beta}-\Gamma(\gamma)^\delta_{\alpha\beta}u_\delta\bigr)
   +\ell\frac{\psi'}{\psi}u_r\delta_{ab}+u\delta_{ab},\\
 W_{an}
 &=\frac{\ell}{\psi}E_a^\alpha u_{\alpha r}
   +\frac1\psi E_a^\alpha
     \left(\ell'-\ell\frac{\psi'}{\psi}\right)u_\alpha,\\
 W_{nn}&=\ell u_{rr}+2\ell'u_r+u.
\end{align*}
Thus, in ordinary coordinate Hessians,
\[
 W=P D^2u P^T+L,
 \qquad
 P=\sqrt\ell\begin{pmatrix}\psi^{-1}E&0\\0&1\end{pmatrix},
\]
where $E=(E_a^\alpha)$. The displayed formulas specify every entry of
$L$. The matrix $P$ is uniformly invertible and smooth on each fixed
coordinate collar. Because $L$ depends smoothly on the fixed geometry,
$u$, and $Du$, its Lipschitz norm on the half-collar is bounded by the
already available $C^2$ norm of $u$.

The needed boundary cancellation is exact:
\[
 \ell(0)=\sin^2\theta,\quad
 \ell'(0)=-\cos\theta\sin\theta,\quad
 \frac{\psi'(0)}{\psi(0)}=-\cot\theta,
 \quad
 \ell'(0)-\ell(0)\frac{\psi'(0)}{\psi(0)}=0.
\]
Consequently $L_{an}(y,0)=0$. Also $u_r(y,0)=0$, and differentiation
along the boundary gives $u_{\alpha r}(y,0)=0$.

Write $R=\operatorname{diag}(1,\ldots,1,-1)$. Extend $u$ evenly by
$\widetilde u(y,r)=u(y,|r|)$. The last two identities imply that this
extension is $C^2$ across $r=0$. For $r<0$ define
\[
 \widetilde P(y,r)=P(y,-r),\qquad
 \widetilde L(y,r)=R L(y,-r)R,
 \qquad \widetilde g(y,r)=g(y,-r),
\]
and keep the original values for $r\geq0$. Since $P$ is block diagonal,
it commutes with $R$. The tangential and normal diagonal blocks of $L$
are reflected evenly, while its mixed block is reflected oddly. Since
$L_{an}(y,0)=0$, $\widetilde L$ is Lipschitz across the interface, with
its Lipschitz constant bounded by the same fixed data. The even
extension $\widetilde P$ is also Lipschitz and uniformly invertible,
and $\widetilde g$ is uniformly $C^\alpha$. For $r<0$,
\[
 \widetilde P D^2\widetilde u\widetilde P^T+\widetilde L
 =R W(y,-r)R.
\]
Thus this reflected matrix has eigenvalues in $[c,C]$ everywhere.

To put the operator into a globally uniformly elliptic class, let
$F(A)=\bar\sigma_k(A)^{1/k}$ and
$\mathcal K=\{A:(c/2)I\leq A\leq2CI\}$. Define
\[
 \widehat F(A)
 =\inf_{B\in\mathcal K}
   \bigl\{F(B)+DF(B):(A-B)\bigr\},
 \qquad A\in\operatorname{Sym}(n).
\]
Concavity of $F$ on the positive cone shows that
$\widehat F=F$ on $\mathcal K$. Since $DF(B)$ has positive eigenvalues
bounded above and below uniformly for $B\in\mathcal K$, there are
$0<\lambda_0\leq\Lambda_0$ such that, for every $N\geq0$,
\[
 \lambda_0\operatorname{tr}N
 \leq\widehat F(A+N)-\widehat F(A)
 \leq\Lambda_0\operatorname{tr}N.
\]
The operator $\widehat F$ is concave, globally Lipschitz, and invariant
under orthogonal conjugation: this last fact follows because both
$F$ and $\mathcal K$ are invariant under such conjugation.

Now set
\[
 \mathfrak F(x,M)
 =\widehat F\bigl(\widetilde P(x)M\widetilde P(x)^T
                      +\widetilde L(x)\bigr).
\]
Uniform invertibility of $\widetilde P$ gives uniform ellipticity of
$\mathfrak F$; concavity is preserved under this affine substitution.
Its Lipschitz dependence on $x$ satisfies
\[
 |\mathfrak F(x,M)-\mathfrak F(y,M)|
 \leq C_0(1+\|M\|)|x-y|,
\]
where $C_0$ depends only on the existing $C^2$ bound and the fixed data.
Moreover $|\mathfrak F(x,0)|$ is bounded by the same data. Orthogonal
invariance and agreement with $F$ on $\mathcal K$ show that
\[
 \mathfrak F(x,D^2\widetilde u)=\widetilde g(x)
\]
on the doubled collar; at $r=0$ the equality follows by continuity.
Application of Lemma~\ref{lem-EKH} on balls in this doubled collar therefore yields
a uniform $C^{2,\beta}$ estimate up to the original boundary.
The same construction without reflection gives the interior estimates.
A finite fixed coordinate covering, followed by multiplication by
$\ell$, proves the asserted estimate for $h$.
\end{proof}

The proof uses only the existing $C^2$ bound to control the Lipschitz
coefficients of the reflected operator. Thus the reduction to Lemma~\ref{lem-EKH}
does not presuppose the second-derivative H\"older estimate that it
seeks to establish.

\begin{proof}[Proof of the long-time existence part of Theorem~\ref{thm-main}]
Let $[0,T_{\max})$ be the maximal smooth strictly convex existence interval. By
Proposition~\ref{prop-C2},
\[
c\sig\le W[h]\le C\sig
\qquad\hbox{on }\Ct\times[0,T_{\max}),
\]
and by Proposition~\ref{prop-C0},
\[
0<c\le h\le C.
\]
Since $\ell/h$ stays in a compact subset of $(0,\infty)$ and $f,\varphi$ are fixed positive smooth
data, the mobility $\Theta$ also satisfies
\[
0<c_\Theta\le\Theta\le C_\Theta.
\]
On the compact matrix set
\[
\{A:\ cI\le A\le CI\},
\]
the derivative matrix $S^{ij}(A)$ is positive definite with eigenvalues bounded above and below.
Consequently $\Theta S^{ij}$ is uniformly elliptic on the entire interval
$[0,T_{\max})$.

The boundary operator $\nabla_\mu h-\cot\theta\,h$ has gradient derivative $\mu$, so it is
uniformly oblique; when $\theta=\pi/2$ it is the homogeneous Neumann condition.

We now explain the regularity bootstrap. The speed ratio solves the linear uniformly parabolic
equation with drift
\[
Q_t-\Theta S^{ij}Q_{ij}
-
2\Theta S^{ij}(\log h)_iQ_j
=
(k+1-a_\varphi(\ell/h))Q(Q-1),
\qquad Q_\mu=0.
\]
The coefficients are bounded on every positive-time strip and the right-hand side is bounded.
The first assertion of Lemma~\ref{lem-parabolic}, applied with
homogeneous Neumann data, therefore gives, for each $\tau>0$,
\[
\|Q\|_{C^{\alpha,\alpha/2}(\overline\Ct\times[\tau,T_{\max}))}\le C_\tau.
\]
At a fixed time the equation
\[
\bsig_k(W[h])^{1/k}
=
\left(
\frac{fQ}{h\varphi(\ell/h)}
\right)^{1/k}
\]
is used in the regularity argument below.

\paragraph{Application to the flow.}
On each fixed positive-time strip, the previously established estimates
give $\|h\|_{C^2}\leq M_0$, $c\sigma\leq W[h]\leq C\sigma$, and
\[
 g(x,t)=\left(\frac{f(x)Q(x,t)}{h(x,t)\varphi(\ell(x)/h(x,t))}\right)^{1/k}.
\]
The bounds for $h$ and $Q$, including their positive lower bounds, together with
the spatial $C^1$ bound for $h$ and the spatial $C^\alpha$ bound for
$Q$, give $\sup_t\|g(\cdot,t)\|_{C^\alpha}\leq M_1$. 

The boundary second-derivative estimate therefore gives a uniform
spatial $C^{2,\gamma}$ bound for some $0<\gamma<\alpha$, where
$\alpha$ denotes the preceding parabolic H\"older exponent for $Q$.
No spatial derivative of $Q$ is required in this step.

Since $h_t=h(Q-1)$, we have $\|h_t\|_\infty\leq C$.
Fix $0<\tau<T_{\max}$ and use the spatial estimate on
$[\tau/2,T_{\max})$. For $s,t$ in this interval, interpolation
on the fixed smooth cap gives
\[
 \|D^2h(t)-D^2h(s)\|_\infty
 \leq C\|h(t)-h(s)\|_\infty^{\gamma/(2+\gamma)}
       \|h(t)-h(s)\|_{C^{2,\gamma}}^{2/(2+\gamma)}
 \leq C_\tau|t-s|^{\gamma/(2+\gamma)}.
\]
Choose
\[
 0<\beta\leq\min\{\alpha,2\gamma/(2+\gamma)\}.
\]
The coefficients $\Theta S^{ij}$ and
$2\Theta S^{ij}(\log h)_i$, and the source in the $Q$-equation,
are uniformly $C^{\beta,\beta/2}$ on positive-time strips.
Lemma~\ref{lem-parabolic}(2) with homogeneous Neumann data then
gives uniform $C^{2+\beta,1+\beta/2}$ estimates for $Q$ on
smaller strips. All these estimates are applied on backward strips
of a fixed length depending on $\tau$, with nested earlier strips
contained in $[\tau/2,T_{\max})$. Their constants are independent
of the strip's location and of $T_{\max}$. The already smooth local
solution covers the initial layer.

Also $h_t=h(Q-1)$ is uniformly $C^{\beta,\beta/2}$.
Together with the spatial and time bounds for $D^2h$, this yields
$h\in C^{2+\beta,1+\beta/2}$ on the smaller strips.
We estimate the time modulus of $Dh$ using centered spatial differences, which have
an error $O(r^{1+\beta})$, while their change between times $s,t$
is bounded by $C|t-s|r^{\beta-1}$ using the spatial
$C^\beta$ bound for $h_t$. Taking $r=|t-s|^{1/2}$ gives the
required $|t-s|^{(1+\beta)/2}$ bound. A fixed smooth-domain
extension operator gives the same estimate in boundary charts.

For clarity, the boundary bootstrap can be carried out with time and
tangential derivatives only. Put \(w=h_t\). Then
\[
 w_t=\Theta S^{ij}\nabla_i\nabla_j w+
 \bigl(\Theta\operatorname{tr}_{\sigma}S
          +(2-a_\varphi(\ell/h))Q-1\bigr)w,
 \qquad w_\mu=\cot\theta\,w.
\]
Its coefficients are uniformly H\"older, and \(w\) is uniformly bounded;
linear parabolic oblique Schauder estimates apply on a smaller strip.
In boundary Fermi coordinates \((y,r)\), with \(r\) inward distance,
the boundary condition is \(h_r+\cot\theta\,h=0\). Consequently each
tangential coordinate derivative \(v_\alpha=\partial_{y^\alpha}h\)
satisfies the same homogeneous Robin condition. Differentiating the
interior equation yields a linear uniformly parabolic equation for
\(v_\alpha\), whose source consists only of derivatives of the coefficients determined
by the fixed geometry and data, together with derivatives of \(h\) of order at most
two; these are uniformly H\"older at the first step. Thus the same
linear estimate controls all third spatial derivatives having a
tangential index. Uniform ellipticity allows the original equation
to be solved locally for \(h_{rr}\); differentiating this identity,
using the estimates for \(h_t\), recovers the missing pure normal
derivatives. Repeat this argument with time and tangential derivatives,
on nested strips, to obtain all orders. At each fixed order the
coefficients and differentiated boundary data required by the linear
Schauder estimates are controlled by the preceding steps. Arbitrary normal
derivatives are never assigned the original Robin boundary condition.
This yields \(\sup_{\tau\leq t<T_{\max}}\|h(t)\|_{C^m}\leq C_{m,\tau}\).
If \(T_{\max}<\infty\), the same bounds and the equation give uniform
\(C^m\) bounds for \(h_t\). Thus \(h(t)\) is Cauchy in every
\(C^m\) as \(t\uparrow T_{\max}\). Repeatedly differentiating the equation
also passes every compatibility identity to the limiting jet.

Under this supposition, the limiting support function is smooth,
positive, and strictly convex by the uniform bounds. The local
existence assertion in Lemma~\ref{lem-parabolic} extends the solution through
$T_{\max}$, contradicting maximality. Hence $T_{\max}=\infty$.

\end{proof}

\subsection{Uniqueness of the stationary solution}\label{subsec-uniqueness}

The supercritical growth condition also gives a simple comparison proof of uniqueness.

\begin{proposition}[Stationary uniqueness]\label{prop-stationary-unique}
Assume
\[
a_\varphi(s)\ge k+1+\delta
\qquad(s>0).
\]
Then \eqref{stationary-eqn} has at most one positive strictly convex solution.
For $k=1$, the same conclusion holds for positive solutions satisfying
$\bsig_1(W[h])>0$.
\end{proposition}

\begin{proof}
Let $h_1,h_2$ be two solutions and set
\[
c=\max_{\overline\Ct}\frac{h_1}{h_2}.
\]
Suppose $c>1$ and let $x_0$ be a maximum point. Then
$w=h_1-ch_2$ has maximum value zero at $x_0$. At an interior maximum,
the usual second-derivative test gives $Dw=0$ and $D^2w\le0$.
% R12 | Boundary maximum test used in uniqueness
% Status: verified optional completion
At a boundary maximum of \(w=h_1-ch_2\),
\(w=0\) and the Robin condition gives \(\nabla_\mu w=0\).
All tangential first derivatives also vanish, so \(Dw=0\).
Taylor expansion along every inward-pointing direction gives
\(D^2w(v,v)\leq0\). Since a quadratic form is even, the same
inequality holds for outward-pointing directions; the inequality in tangential
directions follows by continuity. Hence \(D^2w\leq0\) on the whole
tangent space and \(W[h_1]\leq cW[h_2]\) at the maximum point.

\[
W[h_1](x_0)\le c\,W[h_2](x_0).
\]
For strictly convex matrices, monotonicity and $k$-homogeneity imply
\[
\bsig_k(W[h_1](x_0))
\le
c^k\bsig_k(W[h_2](x_0)).
\]
For $k=1$ the same conclusion follows simply by taking traces.

Since $h_1(x_0)=ch_2(x_0)$, the two stationary equations give
\begin{align*}
1
&=
\frac{
h_1\varphi(\ell/h_1)\bsig_k(W[h_1])
}{
h_2\varphi(\ell/h_2)\bsig_k(W[h_2])
}(x_0)\\
&\le
c^{k+1}
\frac{
\varphi\!\left(\ell(x_0)/(c h_2(x_0))\right)
}{
\varphi\!\left(\ell(x_0)/h_2(x_0)\right)
}.
\end{align*}
The supercritical hypothesis implies
\[
\log
\frac{
\varphi\!\left(\ell(x_0)/h_2(x_0)\right)
}{
\varphi\!\left(\ell(x_0)/(c h_2(x_0))\right)
}
=
\int_{\ell(x_0)/(c h_2(x_0))}^{\ell(x_0)/h_2(x_0)}
\frac{a_\varphi(s)}s\,\dd s
\ge
(k+1+\delta)\log c.
\]
Hence
\[
\frac{
\varphi\!\left(\ell(x_0)/(c h_2(x_0))\right)
}{
\varphi\!\left(\ell(x_0)/h_2(x_0)\right)
}
\le
c^{-(k+1+\delta)}
\]
and consequently
\[
1\le c^{-\delta}<1,
\]
a contradiction. Thus $h_1\le h_2$. Interchanging the two solutions gives
$h_2\le h_1$, and therefore $h_1=h_2$.
\end{proof}

\subsection{Exponential convergence}\label{subsec-convergence}

Let
\[
M(t)=\max_{\overline\Ct}Q(\cdot,t),
\qquad
m(t)=\min_{\overline\Ct}Q(\cdot,t).
\]
From \eqref{Q-evolution}, if $M(t)>1$,
\[
D^+(M-1)
\le
-\delta M(M-1)
\le-\delta(M-1).
\]
If $m(t)<1$, Proposition~\ref{prop-Q} gives $m(t)\ge c_Q>0$ and
\[
D^+(1-m)
\le
-\delta c_Q(1-m).
\]
Hence there are $C,c>0$ such that
\begin{equation}\label{Q-exp}
\|Q(\cdot,t)-1\|_{L^\infty(\Ct)}
\le Ce^{-ct}.
\end{equation}

Since
\[
\partial_t\log h=Q-1,
\]
\eqref{Q-exp} implies
\[
\int_0^\infty
\|\partial_t\log h(\cdot,t)\|_{L^\infty}\dd t<\infty.
\]
Therefore $h(\cdot,t)$ converges uniformly to a positive function $h_\infty$.
The time-uniform higher-order estimates from Subsection~\ref{subsec-longtime} imply precompactness in every
$C^m$, and the uniqueness of the $C^0$ limit upgrades the convergence to
\[
h(\cdot,t)\to h_\infty
\quad\hbox{in }C^\infty(\overline\Ct).
\]
Passing to the limit in $Q\to1$ gives
\[
h_\infty\varphi(\ell/h_\infty)
\bsig_k(W[h_\infty])=f.
\]
The Robin condition passes to the limit as well. Proposition~\ref{prop-stationary-unique}
shows that the limiting stationary solution is unique. This proves the convergence part of
Theorem~\ref{thm-main}.

\subsection{The all-angle scalar flow when $k=1$}\label{subsec-k1}

We now prove Theorem~\ref{thm-k1-full-angle}. No convexity estimate is used in this section.

For $k=1$,
\[
\bsig_1(W[h])=\frac1n(\Delta h+nh).
\]
Set $u=h/\ell$. By \eqref{Btheta},
\[
\bsig_1(W[h])
=
\frac1n\left(
\ell\Delta u
+
2\langle\nabla\ell,\nabla u\rangle
+
nu
\right).
\]
Thus the flow becomes
\begin{equation}\label{k1-u-flow}
\begin{aligned}
u_t={}&
a(x,u)\Delta u+b(x,u)\cdot\nabla u+c(x,u),\\
u_\mu={}&0,
\end{aligned}
\end{equation}
where
\begin{align}
a(x,u)
&=
\frac{\ell^2u^2\varphi(1/u)}{nf},
\label{k1-a}\\
b(x,u)
&=
\frac{2\ell u^2\varphi(1/u)}{nf}\,\nabla\ell,
\notag\\
c(x,u)
&=
\frac{\ell u^3\varphi(1/u)}f-u.
\notag
\end{align}
For every fixed $0<\theta<\pi$, \eqref{ell-full-angle-bounds} gives $\ell_->0$.

\subsubsection{Uniform scalar bounds and global parabolicity}

Define
\[
g_1(s)=s^2\varphi(1/s).
\]
Then
\[
\frac{g_1'(s)}{g_1(s)}
=
\frac{2-a_\varphi(1/s)}s
\le-\frac\delta s.
\]
Hence $g_1(s)\to\infty$ as $s\downarrow0$ and $g_1(s)\to0$ as $s\uparrow\infty$.

At a spatial maximum of $u$, the drift term in \eqref{k1-u-flow} vanishes and $\Delta u\le0$.
Therefore
\[
D^+u_{\max}
\le
u_{\max}
\left[
\frac{\ell_{\max}}{f_{\min}}g_1(u_{\max})-1
\right].
\]
At a spatial minimum,
\[
D_+u_{\min}
\ge
u_{\min}
\left[
\frac{\ell_{\min}}{f_{\max}}g_1(u_{\min})-1
\right].
\]
The same argument applies at the boundary because $u_\mu=0$. Consequently
\begin{equation}\label{k1-C0}
0<c\le u\le C
\qquad\text{for all times of existence}.
\end{equation}

By \eqref{k1-a}, \eqref{ell-full-angle-bounds}, and \eqref{k1-C0},
\[
0<a_0\le a(x,u)\le a_1<\infty.
\]
% R11 | k=1 all-angle scalar regularity chain
We now give the regularity argument for the scalar equation.
Along the given solution, it can be written as the linear equation
\[
 u_t-a(x,u)\Delta u-b(x,u)\cdot\nabla u=c(x,u),\qquad u_\mu=0.
\]
The bounds \(c_0\leq u\leq C_0\) imply uniform ellipticity of
\(a(x,u)\), and uniform bounds for \(a(x,u),b(x,u),c(x,u)\).
At this stage, derivatives of these composed coefficients are not assumed
bounded: only the derivatives of the functions \(a,b,c\) in
\((x,z)\in\overline{\mathcal C_\theta}\times[c_0,C_0]\) are bounded.
In Fermi coordinates near the boundary, extend \(u\) evenly across
\(x_n=0\). Since \(u_{x_n}=0\), this is a classical \(C^2\)
extension. Reflect the principal matrix by \(A\mapsto RAR\),
the drift by \(b\mapsto Rb\), and the source evenly, where
\(R=\operatorname{diag}(1,\ldots,1,-1)\).
Uniform ellipticity and all required coefficient bounds are preserved.
Interior parabolic Krylov--Safonov estimates for bounded measurable
principal coefficients, bounded drift, and bounded source now yield a
uniform \(C^{\alpha,\alpha/2}\) bound for \(u\) on positive-time
strips. No preliminary gradient estimate is required.
Consequently the composed coefficients are uniformly
\(C^{\alpha,\alpha/2}\). The linear parabolic Schauder estimate
with smooth homogeneous Neumann boundary data yields a uniform
\(C^{2+\alpha,1+\alpha/2}\) bound on a smaller strip.
Higher-order Neumann Schauder estimates, with the smooth dependence of
\(a,b,c\) on \((x,u)\), give bounds for all spatial and time derivatives by
induction. For each fixed \(\tau>0\), use nested backward strips of
length bounded below by a positive number depending on \(\tau\)
(and on the desired derivative order); this makes the constants
independent of the position of the strip and of the maximal time.
Thus, on every maximal smooth interval \([0,T_{\max})\),
\begin{equation}\label{k1-high}
 \sup_{t\in[\tau,T_{\max})}\|u(t)\|_{C^m}\leq C_{m,\tau}.
\end{equation}
If \(T_{\max}<\infty\), these estimates and the equation produce a smooth
positive limiting jet at \(T_{\max}\), preserving the compatibility hierarchy,
and the short-time theorem extends the solution. This proves global
existence. All constants may depend on the fixed angle; no uniformity
as \(\theta\downarrow0\) or \(\theta\uparrow\pi\) is asserted.

\subsubsection{The speed ratio and convergence}

For $k=1$, Lemma~\ref{lem-Q} becomes
\begin{equation*}
Q_t-a(x,h)\Delta Q
=
2a(x,h)\langle\nabla\log h,\nabla Q\rangle
+
\bigl(2-a_\varphi(\ell/h)\bigr)Q(Q-1),
\end{equation*}
where
\[
a(x,h)=\frac{\Theta}{n}>0,
\qquad
Q_\mu=0.
\]
The initial assumption $\bsig_1(W[h_0])>0$ is exactly $Q(\cdot,0)>0$.
The maximum principle therefore gives
\[
0<c_Q\le Q\le C_Q.
\]
Since $2-a_\varphi\le-\delta$, the same logistic comparison as in Subsection~\ref{subsec-convergence} gives
\[
\|Q(\cdot,t)-1\|_{L^\infty}\le Ce^{-ct}.
\]
Because $\partial_t\log h=Q-1$, $\log h(\cdot,t)$ is Cauchy in $C^0$ and converges
exponentially. The uniform estimates \eqref{k1-high} and interpolation upgrade the convergence to
$C^\infty$. Passing to the limit in the equation gives the stationary $k=1$ equation.

Finally, Proposition~\ref{prop-stationary-unique} gives uniqueness of the positive stationary
solution with $\bsig_1(W)>0$. This proves Theorem~\ref{thm-k1-full-angle}. Notice that nowhere in
the scalar continuation argument did we use the sign of $\cot\theta$. In the obtuse range, however, this scalar continuation argument does not establish
the separate geometric assertion that the full tensor $W[h]$ stays positive definite.

\subsection{The natural boundary condition and the contact-angle problem}\label{subsec-angle}

The quantitative hypothesis \eqref{boundary-quant} is a convenient sufficient condition for a
time-uniform non-collapsing estimate, but it is stronger than the natural static condition
\[
H_\mu\ge0.
\]
Under this weaker condition, Theorem~\ref{thm-preserve} shows that strict convexity cannot be lost
% R13 | Remove undefined eta and A
% Status: verified minor notation correction
at any finite smooth time in the non-obtuse range.
For a normal minimum-radius support $v=hW(\mu,\mu)$, the formulas in Lemma~\ref{lem-q-boundary}, with $k=1$ treated separately, give
\[
v_\mu+(n-1)\cot\theta\,v\ge0.
\]

Thus the remaining issue under the natural boundary sign is quantitative rather than qualitative:
one seeks a time-independent positive lower bound for $q$ under this Robin sign alone.

The angle restriction is independent of this optimization. At the right angle, $\ell\equiv1$ and
$\cot\theta=0$, so all boundary signs used above become Neumann signs and the global theorem
remains valid. In contrast, for $\pi/2<\theta<\pi$, a tangential minimum-radius support satisfies
\[
v_\mu=h\cot\theta\,\lambda_n<0.
\]
This has the wrong sign for the boundary barriers used in the finite-time preservation and
quantitative non-collapsing arguments. Hence neither proof closes by the present boundary
mechanism. The obtuse-angle problem therefore requires a different boundary curvature
estimate.

For $k=1$ the scalar equation is uniformly parabolic after the $C^0$ estimate, and
Theorem~\ref{thm-k1-full-angle} gives an all-angle scalar flow
for $0<\theta<\pi$.
What remains open in the obtuse range is the geometric strengthening: preservation of the full
matrix inequality $W>0$ and hence an interpretation as a strictly convex capillary hypersurface
for all times.

\begin{remark}
Theorems~\ref{thm-main} and \ref{thm-preserve} are stated separately. The global convergence
theorem uses a boundary condition for which the non-collapsing argument closes directly,
whereas the optimal natural boundary hypothesis and the obtuse-angle extension remain
separate technical problems.
\end{remark}

\medskip
\noindent\textbf{Acknowledgments:}
Hua Zhu was supported by the National Natural Science Foundation of
China (Grant No. 12501273, 12661041) and the Research Foundation
of Southwest University of Science and Technology
(Grant No. 25zx7153, 25zx7154). Peihe Wang is supported by National Natural Science Foundation of Shandong Province(ZR2026MS0039).
\medskip

\section*{Declaration on generative AI assistance}

During the preparation of this manuscript, generative AI tools (GPT) were used to assist with language editing, literature searches, and the checking of mathematical arguments. The authors are responsible for the correctness of the mathematical arguments.


\begin{thebibliography}{99}

\bibitem{BISOrlicz}
P. Bryan, M. N. Ivaki, and J. Scheuer,
\emph{Orlicz--Minkowski flows},
Calc. Var. Partial Differential Equations \textbf{60} (2021), Article 41.


\bibitem{BISParabolic}
P. Bryan, M. N. Ivaki, and J. Scheuer,
\emph{Parabolic approaches to curvature equations},
Nonlinear Anal. \textbf{203} (2021), Article 112174.


\bibitem{BISConstantRank}
P. Bryan, M. N. Ivaki, and J. Scheuer,
\emph{Constant rank theorems for curvature problems via a viscosity approach},
Calc. Var. Partial Differential Equations \textbf{62} (2023), Article 98.


\bibitem{CNS}
L. Caffarelli, L. Nirenberg, and J. Spruck,
\emph{The Dirichlet problem for nonlinear second order elliptic equations, III:
Functions of the eigenvalues of the Hessian},
Acta Math. \textbf{155} (1985), 261--301.


\bibitem{GardnerHugWeil}
R. J. Gardner, D. Hug, and W. Weil,
\emph{The Orlicz--Brunn--Minkowski theory: A general framework, additions, and inequalities},
J. Differential Geom. \textbf{97} (2014), 427--476.


\bibitem{GuanMa}
P. Guan and X.-N. Ma,
\emph{The Christoffel--Minkowski problem I: Convexity of solutions of a Hessian equation},
Invent. Math. \textbf{151} (2003), 553--577.


\bibitem{GuanMaZhou}
P. Guan, X.-N. Ma, and F. Zhou,
\emph{The Christoffel--Minkowski problem III: Existence and convexity of admissible solutions},
Comm. Pure Appl. Math. \textbf{59} (2006), 1352--1376.


\bibitem{HaberlLutwakYangZhang}
C. Haberl, E. Lutwak, D. Yang, and G. Zhang,
\emph{The even Orlicz--Minkowski problem},
Adv. Math. \textbf{224} (2010), 2485--2510.


\bibitem{HuIvaki}
Y. Hu and M. N. Ivaki,
\emph{Capillary $L_p$-Christoffel--Minkowski problem},
\href{https://arxiv.org/abs/2512.15464v2}{arXiv:2512.15464v2}, 2025.

\bibitem{HuIvakiCurvature}
Y. Hu and M. N. Ivaki,
\emph{Capillary $L_p$-curvature problem},
Adv. Calc. Var. \textbf{19} (2026), 699--705.


\bibitem{HISCapillaryCM}
Y. Hu, M. N. Ivaki, and J. Scheuer,
\emph{Capillary Christoffel--Minkowski problem},
\href{https://arxiv.org/abs/2504.09320v1}{arXiv:2504.09320v1}, 2025.

\bibitem{Ivaki2019}
M. N. Ivaki,
\emph{Deforming a hypersurface by principal radii of curvature and support function},
Calc. Var. Partial Differential Equations \textbf{58} (2019), Article 1.


\bibitem{KrylovSafonov}
N. V. Krylov and M. V. Safonov,
\emph{A certain property of solutions of parabolic equations with measurable coefficients},
Math. USSR-Izv. \textbf{16} (1981), 151--164.


\bibitem{LSU}
O. A. Ladyzhenskaya, V. A. Solonnikov, and N. N. Ural'tseva,
\emph{Linear and Quasi-linear Equations of Parabolic Type},
Translations of Mathematical Monographs, vol.~23,
American Mathematical Society, Providence, RI, 1968.


\bibitem{LiLiuFlow}
G. Li and C. Liu,
\emph{Capillary Orlicz--Minkowski flow in the upper half-space},
\href{https://arxiv.org/abs/2601.14659v1}{arXiv:2601.14659v1}, 2026.

\bibitem{LiebermanTrudinger}
G. M. Lieberman and N. S. Trudinger,
\emph{Nonlinear oblique boundary value problems for nonlinear elliptic equations},
Trans. Amer. Math. Soc. \textbf{295} (1986), 509--546.


\bibitem{Lutwak93}
E. Lutwak,
\emph{The Brunn--Minkowski--Firey theory. I. Mixed volumes and the Minkowski problem},
J. Differential Geom. \textbf{38} (1993), 131--150.


\bibitem{MWWMinkowski}
X. Mei, G. Wang, and L. Weng,
\emph{The capillary Minkowski problem},
Adv. Math. \textbf{469} (2025), Article 110230.


\bibitem{MWWLpCurvature}
X. Mei, G. Wang, and L. Weng,
\emph{Prescribed $L_p$ curvature problem for convex capillary hypersurface},
\href{https://arxiv.org/abs/2512.16686}{arXiv:2512.16686}, 2025.

\bibitem{MWWChristoffel}
X. Mei, G. Wang, and L. Weng,
\emph{The capillary Christoffel--Minkowski problem},
Calc. Var. Partial Differential Equations \textbf{65} (2026), Article 186.


\bibitem{MWWLpMinkowski}
X. Mei, G. Wang, and L. Weng,
\emph{The capillary $L_p$-Minkowski problem},
Trans. Amer. Math. Soc., to appear;
\href{https://arxiv.org/abs/2505.07746v3}{arXiv:2505.07746v3}.


\bibitem{SafonovBVP}
M. V. Safonov,
\emph{On the boundary value problems for fully nonlinear elliptic equations of second order},
author manuscript, 1994.

\bibitem{WWXAF}
G. Wang, L. Weng, and C. Xia,
\emph{Alexandrov--Fenchel inequalities for convex hypersurfaces in the half-space with capillary boundary},
Math. Ann. \textbf{388} (2024), 2121--2154.


\bibitem{WWXMinkowski}
G. Wang, L. Weng, and C. Xia,
\emph{A Minkowski-type inequality for capillary hypersurfaces in a half-space},
J. Funct. Anal. \textbf{287} (2024), Article 110496.


\bibitem{WangZhu}
X. Wang and B. Zhu,
\emph{The capillary Orlicz--Minkowski problem},
J. Funct. Anal. \textbf{291} (2026), Article 111574.


\bibitem{XiJinLeng}
D. Xi, H. Jin, and G. Leng,
\emph{The Orlicz Brunn--Minkowski inequality},
Adv. Math. \textbf{260} (2014), 350--374.


\bibitem{ZhuStatic}
H. Zhu,
\emph{The capillary Orlicz--Christoffel--Minkowski problem},
\href{https://arxiv.org/abs/2611.18299v1}{arXiv:2611.18299v2}, 2026.

\end{thebibliography}
\end{document}